\documentclass[letterpaper,11pt]{elsarticle}
\usepackage{amsmath}
\usepackage{times}
\usepackage{hyperref}
\usepackage{amssymb}
\usepackage{graphicx}
\usepackage{subcaption}
\usepackage{multirow}
\usepackage{booktabs}
\usepackage{listings}

\DeclareFixedFont{\Tifont}{T1}{ptm}{b}{it}{16pt}
\usepackage{setspace}
\usepackage{amssymb}
\usepackage{algorithm}
\usepackage{algorithmic}

\newcommand{\expnumber}[2]{{#1}\mathrm{e}{#2}}

\newcommand{\figref}[1]{\figurename~\ref{#1}}

\newcommand{\vf}{{\mathbf{f}}}

\newcommand{\vn}{{\mathbf{n}}}

\newcommand{\pddt}{{\frac{\partial}{\partial t}}}

\usepackage{xcolor}
\definecolor{deepblue}{rgb}{0,0,0.5}
\newcommand\cppstyle{\lstset{
language=C++,
basicstyle=\footnotesize\ttfamily,
breaklines=true,
keywordstyle=\bfseries\color{deepblue},
}}
\newcommand\pythonstyle{\lstset{
language=Python,
basicstyle=\footnotesize\ttfamily,
breaklines=true,
keywordstyle=\bfseries\color{deepblue},
}}

\lstnewenvironment{cppcode}[1][]{\cppstyle \lstset{#1}}{}
\lstnewenvironment{pythoncode}[1][]{\pythonstyle \lstset{#1}}{}

\newcommand\cppinline[1]{{\cppstyle\lstinline!#1!}}
\newcommand\pythoninline[1]{{\pythonstyle\lstinline!#1!}}

\begin{document}

\title{Domain-specific implementation of high-order \\
Discontinuous Galerkin methods in spherical geometry}

\author[1]{Kalman Szenes}%
\ead{kalman.szenes@math.ethz.ch}
\author[2]{Niccol{\`o} Discacciati}
\ead{niccolo.discacciati@epfl.ch}
\author[3]{Luca Bonaventura}
\ead{luca.bonaventura@polimi.it}
\author[4]{William Sawyer\corref{cor4}}%
\ead{william.sawyer@cscs.ch}
\cortext[cor4]{Corresponding author}

\affiliation[1]{organization={Swiss Federal Institute of Technology Zurich},
addressline={Raemistrasse 101},
postcode={8092},
city={Zurich},
country={Switzerland}}
\affiliation[2]{organization={Swiss Federal Institute of Technology Lausanne},
addressline={Route Cantonale},
postcode={1015},
city={Lausanne},
country={Switzerland}}
\affiliation[3]{organization={Politecnico di Milano},
addressline={Piazza Leonardo da Vinci 32},
postcode={20133},
city={Milan},
country={Italy}}
\affiliation[4]{organization={Swiss National Supercomputing Centre},
addressline={Via Trevano 131},
postcode={6900},
city={Lugano},
country={Switzerland}}


\begin{abstract}
We assess two domain-specific languages included in the GridTools ecosystem as tools for implementing a high-order Discontinuous Galerkin discretization of the shallow water equations. Equations in spherical geometry are considered, thus providing a blueprint for the application of domain-specific languages to the development of global atmospheric models. The  results  demonstrate that  domain-specific languages designed for finite difference/volume methods can be successfully extended to implement a Discontinuous Galerkin solver.

\end{abstract}

  \maketitle


\section{Introduction}  \label{sec:intro}


It has always been challenging for numerical mathematicians to implement new algorithms in a way that can attain the best possible performance of the underlying platform.  The task has become more difficult with the evolution of new computing architectures, such as Graphics Processing Units (GPUs) or Field-programmable Gate Arrays (FPGAs).  Scientific programmers are forced to learn computing concepts well outside of their original domain.

Domain-specific languages (DSLs) \cite{fowler2010domain} represent an attempt  to separate the concerns of the domain scientist from the complexities of the underlying computer science issues.  In our case, the designer formulates her problem in terms of numerical operations, including time-stepping algorithms, linear algebra operations, or Partial Differential Equation (PDE) formulations, while a ``backend'' takes care of generating the appropriate code for the architecture.  While programming languages such as C++ or Fortran attempt to abstract away the hardware, they have not kept pace with the emerging hardware complexity and frequently require compiler directives or add-ons to properly exploit the architecture.

Extensive effort has been made in developing frameworks to ease the challenge of solving PDEs in given geometries.  On the far end of the spectrum, software frameworks like FEniCS \cite{alnaes2015fenics} and Firedrake \cite{rathgeber2016firedrake} offer a descriptive language to define the PDE and the given domain, along with initial and boundary conditions.  These frameworks then generate the code to solve the problem using finite element methods. The former gives the user little latitude to test a new numerical technique, making it more appropriate for domain scientists relying on standard, widely-supported methods.  Firedrake, on the other hand, has a different underlying implementation allowing more flexibility to employ code optimizations, various finite element mesh topologies, as well as parallelization features, such as the support of message-passing, multithreading or GPUs. Firedrake utilizes PyOp2~\cite{rathgeber2012PyOP2} -- a full-fledged DSL for the parallel executions of computational kernels on unstructured meshes or graphs -- to allow these interventions. More generally, PyOp2 can be viewed as a DSL for finite element calculations. PyOp2 was subsequently refined into the Psyclone~\cite{Psyclone} code-generation system, which was specifically designed to extend Fortran codes. The latter was then used to implement the LFRic~\cite{kavcic2020lfric} atmospheric model. Psyclone is essentially a source-code generator and  could conceptually address the algorithms we propose here, but we have chosen to evaluate tools which have been developed internally at the Swiss National Supercomputing Centre (CSCS).

Our interests are in the area of weather and climate simulation, where the mesh is often built by extrusion of a two-dimensional horizontal mesh covering, for example, the whole globe.  Our goal is to enable climate and numerical weather prediction (NWP) applications to leverage a variety of architectures by utilizing software backends. A first attempt at a C++-embedded DSL specific to the climate and weather domain was STELLA~\cite{gysi2015stella}, with which the COSMO dynamical core (solver of the non-hydrostatic equations of atmospheric motion) was implemented \cite{thaler2019porting}.  This prototype was completely replaced by GridTools~\cite{afanasyev2021gridtools}, in which the COSMO and NICAM \cite{kunkel2020aimes} models have been implemented. The classic implementation of GridTools assumes a Cartesian grid, which explains the choices  described in Section \ref{sec:math}, however the newest version, released in 2023, also allows for a fully unstructured mesh.

In this paper we evaluate two DSLs included in the GridTools ecosystem, namely the C++-based Galerkin-for-Gridtools (G4GT) and the Python-based GridTools-for-Python (GT4Py)~\cite{GridTools}, as tools for implementing a high-order Discontinuous Galerkin (DG) method for time-dependent problems, see e.g., \cite{giraldo:2020,hesthavenNodalDiscontinuousGalerkin2008}. The GridTools framework was originally designed to support finite difference/volume methods on rectangular grids and, more recently, on unstructured grids. G4GT was a first prototype to extend these tools for finite element problems. GT4Py is a Python layer above GridTools and thus targets the same FD/FV problems, for example \cite{dahm2021gt4py}, but we found that with some extensions it can be reused to implement DG solvers. More specifically, we consider an explicit time discretization and a modal DG spatial discretization for a system of conservation laws.  
Common, but non-trivial, benchmarks, namely linear advection and the shallow water equations, are used to validate our DSL implementations of a DG method.  Extensive literature is available on the results obtained for these benchmarks by a wide range of numerical methods, allowing us to assess the numerical and performance results. One of the two implementations deals with equations in spherical geometry, showing that DSLs can be successfully applied to the development of prototype codes for atmospheric modeling.

The structure of the paper is as follows.
The flux formulation of the shallow water equations in latitude-longitude coordinates follows in Section~\ref{sec:math}.  A preliminary implementation in an early C++-based DSL prototype called G4GT is presented in Section~\ref{sec:impl_G4GT}, while the new, Python-based implementation called GT4Py is discussed at length in Section~\ref{sec:impl_GT4Py}.  
The validation of the resulting implementations is discussed in Section~\ref{sec:validation} and benchmarks are presented in Section~\ref{sec:performance}.  In Section~\ref{sec:conclusions}, we recount our experiences using the DSLs and make suggestions on how to better support finite element codes in these frameworks.

\section{The mathematical model and numerical  discretization approach}  \label{sec:math}

\begin{figure}[htb]
        \centering
        \includegraphics[width=0.2\textwidth]{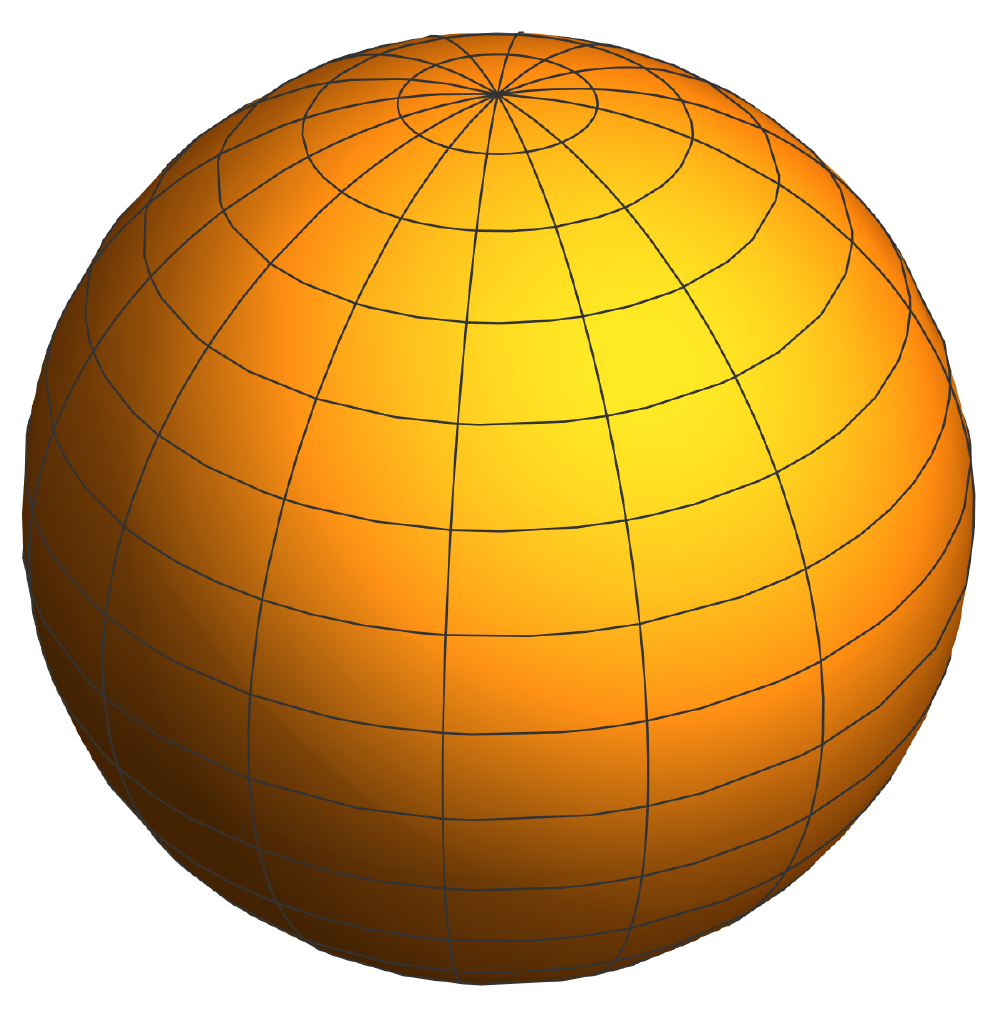}
        \caption[lat-lon]{Latitude-Longitude grid}
        \label{fig:lat-lon}
    \end{figure}

We are concerned with demonstrating the capabilities of DSLs for implementing numerical solutions of conservation laws. 
A system of conservation laws can be written as:
    \begin{equation} 
        \frac{\partial \textbf{u}}{\partial t} + \nabla \cdot \textbf{F}(\textbf{u}) = \textbf{S}(\textbf{u}),  \label{eq:fluxform}
    \end{equation}
\noindent
where $\textbf{u}$ is the vector of conserved variables, $\textbf{F}$ is the flux function and $\textbf{S}$ is the source term.
Equation \eqref{eq:fluxform} becomes well posed once complemented with appropriate initial conditions and boundary conditions, see, for example, the discussion in \cite{levequeFiniteVolumeMethods2002}.
Since our goal is the application of DSL tools to models for weather and climate, we choose as main model equations the shallow water equations (SWEs) on the sphere, which are a common benchmark for numerical models in this area. Various formulations of these equations can be found in \cite{williamsonStandardTestSet1992}.
We consider the Earth's surface as  a sphere of
 radius 
 $R,$ that is parameterized   in latitude - longitude (lat-lon) coordinates as a rectangular domain such that
 the latitude $\theta \in [-\pi/2,\pi/2] $ and the longitude
 $\lambda \in [0,2\pi].$ Denote by $\Omega =7.292\times 10^{-5} \rm \ s^{-1}$   the Earth's rotation rate, by $f=2\Omega \sin \theta $ the Coriolis
 parameter and by  $g=9.81 \ \rm m \ s^{-2}$ the Earth's gravitational acceleration. 
Furthermore, let   $\hat{\boldsymbol \imath }$,
 $ \hat{\boldsymbol \jmath }$   denote the longitudinal and latitudinal  unit vectors, respectively.
  For a generic scalar function $\phi$
 and vector field $\mathbf{w}=w_1\hat{\boldsymbol \imath }+w_2\hat{\boldsymbol \jmath }, $  
   the spherical gradient and divergence are defined as:

 \begin{eqnarray}
 \label{eq:operators}
     \nabla \phi &=& \frac{\hat{\boldsymbol \imath }}{R \cos\theta} \partial_\lambda \phi + \frac{\hat{\boldsymbol \jmath }}{R} \partial_\theta \phi \nonumber\\
     \nabla \cdot \mathbf{w} &=& \cfrac{1}{R \cos\theta} \left[\partial_\lambda(w_1) + \partial_\theta(w_2 \cos\theta)\right].
\end{eqnarray}
 Denoting  then $h$ as the thickness
of a fluid over the spherical surface (assuming flat orography $h_b=0$) and $\mathbf{v}=u\hat{\boldsymbol \imath}+v\hat{\boldsymbol \jmath} $ the velocity field, which is a tangent vector field to the sphere, the shallow water equations
in flux form are written as: 

\begin{eqnarray} 
\label{eq:sph1}
      &&    \partial_t h + \nabla \cdot (h \mathbf{v}) = 0  \nonumber \\
&&          \partial_t (h \mathbf{v}) + \nabla \cdot (h \mathbf{v} \otimes \mathbf{v}) = - f \mathbf{\hat k} \times h \mathbf{v} -  \nabla \left(\frac{gh^2}2\right).
 \end{eqnarray}
It is well known that lat-lon coordinates entail a number of numerical difficulties. Indeed, for a Cartesian lat-lon  mesh such as the one depicted in Figure \ref{fig:lat-lon}, the elements become increasingly distorted as they approach the poles. Moreover, the elements precisely neighboring the poles have a singular edge in physical space (these elements reduce to spherical triangles instead of rectangles).

    \noindent
While several alternatives have been considered in the literature, see, e.g., the review in \cite{bonaventura:2012},
this setting is sufficient for the present purpose of validating  GT4Py, which, as discussed in Section \ref{sec:intro}, could only handle non-Cartesian meshes in the latest release. Furthermore, degree adaptivity techniques as in \cite{tumolo:2015} can help reduce the numerical problems of this extremely simple setting. The equations for the spherical components of the velocity
field can then be derived taking into account the non-inertial
nature of the rotating reference frame on the sphere, see 
again \cite{williamsonStandardTestSet1992} (Equations 6-7). We then obtain:

\begin{eqnarray} 
\label{eq:sph2}
&&      \partial_t (h u) + \nabla \cdot (hu \mathbf{v} ) +\frac{1}{R\cos\theta}\partial_{\lambda}\left(\frac{gh^2}2\right)
= \left(f+\frac uR \tan \theta\right)hv\nonumber \\
&& \partial_t (h v) + \nabla \cdot (hv \mathbf{v} ) +
\frac{1}{R}\partial_{\theta}\left(\frac{gh^2}2\right)
= -\left(f+\frac uR \tan \theta\right)hu.
 \end{eqnarray}
 Since
 \begin{eqnarray}
  \nabla \cdot (hu\mathbf{v}) &=& \cfrac{1}{R \cos\theta} \left[\partial_\lambda(hu^2) + \partial_\theta(huv \cos\theta)\right]\nonumber \\
  \nabla \cdot (hv\mathbf{v}) &=& \cfrac{1}{R \cos\theta} \left[\partial_\lambda(huv) + \partial_\theta(hv^2 \cos\theta)\right]\nonumber
\end{eqnarray}
and also
$$
\frac{\cos \theta}{R}\partial_{\theta}\left(\frac{gh^2}2\right)=\frac{1}{R}\partial_{\theta}\left(\frac{gh^2}{2}\cos \theta\right) + \frac{gh^2}{2 R}\sin\theta,
$$
the SWE can be rewritten component-wise as:
\begin{eqnarray} 
\label{eq:sph3}
             \partial_t (h \cos\theta)  &+& \cfrac{1}{R}\left[ \partial_\lambda \left(h u\right) + \partial_\theta \left(h v \cos\theta\right) \right] = 0 \nonumber \\
         \partial_t (h u \cos\theta)  &+&  \cfrac{1}{R}\left[\partial_\lambda \left(h u^2 + \cfrac{g h^2}{2}\right) + \partial_\theta \left(h u v \cos\theta\right) \right] \nonumber \\
 &=& \left(f\cos\theta+\frac uR \sin \theta  \right)hv \\
         \partial_t (h v \cos\theta) &+& \cfrac{1}{R}\left[ \partial_\lambda \left(h u v\right) + \partial_\theta \left(\left(h v^2 + \cfrac{g h^2}{2}\right) \cos\theta\right) \right] \nonumber \\
         &=&   - \cfrac{g h^2 \sin \theta}{2 R} 
         -\left(f\cos\theta+\frac uR \sin \theta\right)hu.
         \nonumber
    \end{eqnarray}

Periodic boundary conditions are considered in the longitudinal direction, while in the latitudinal direction, the fluxes are set to zero. The SWE in spherical coordinates can therefore be written as the system of conservation laws:
    \begin{equation*}
     \partial_t (\mathbf{U}) + \partial_\lambda (\mathbf{F}(\mathbf{U})) + \partial_\theta (\mathbf{G}(\mathbf{U})) = \mathbf{S}(\mathbf{U}),
    \end{equation*} \label{eq:SWES}
\noindent
where we have defined the conserved quantities, fluxes and sources as:

    \begin{eqnarray}
       \mathbf{U} &=& \begin{pmatrix}
              h \cos \theta \\
              h u \cos \theta \\
              h v \cos \theta
         \end{pmatrix}, 
         \hskip 0.2cm  \mathbf{F(\mathbf{U})} = \cfrac{1}{R} \begin{pmatrix}
               h u \\
               h u^2 + \cfrac{g h^2}{2} \\
               h u v
         \end{pmatrix} \nonumber \\  
         \mathbf{G}(\mathbf{U}) &=& \cfrac{\cos \theta}{R} \begin{pmatrix}
              h v  \\
              h u v \\
              h v^2 + \cfrac{g h^2}{2}
         \end{pmatrix}, \hskip 0.2cm\mathbf{S}(\mathbf{U}) = \begin{pmatrix}
              0 \\
              \left(f\cos\theta+\frac uR \sin \theta  \right)hv \\
               - \cfrac{g h^2 \sin \theta}{2 R} 
         -\left(f\cos\theta+\frac uR \sin \theta\right)hu
          \end{pmatrix}.  
    \end{eqnarray}
        
\noindent   
We then present an overview of the classical DG method chosen for the demonstration of a DSL implementation. A complete description can be found, among many others, in \cite{giraldo:2020,hesthavenNodalDiscontinuousGalerkin2008}.
We consider for simplicity the discretization of a scalar conservation law,

 \begin{equation} \label{eq:cons_law_sc}
        \frac{\partial {u}}{\partial t} + \nabla \cdot \textbf{f}(u) ={s(u)},
    \end{equation}
defined on a two-dimensional rectangular domain. This domain is subdivided into $K$ elements, denoted by $D^k$, for any $k = 1\ldots K$. We restrict our attention to conforming structured meshes composed of rectangular elements $  D^k = [x^k_l, x^k_r] \times [y^k_b , y^k_t ]. $ 
In each element, let $V^k_h$ be the finite-dimensional space of multivariate polynomials up to a given degree $r $ in each spatial dimension:

\begin{equation*}
V^k_h =
  \left\{ v : v = \sum_{i,j=0}^r \alpha_{ij} x^i y^j,
               x \in [x^k_l,x^k_r], y \in [y^k_b , y^k_t]
  \right\}.
\end{equation*}

\noindent
Consequently, the finite-dimensional space in which we seek the solution is the space of discontinuous polynomials defined as
$
V_h = \left\{ 
         v \in L^2( \Omega ) : v|_{D^k} \in V^k_h
      \right\}.
$
The numerical solution can be viewed as the direct sum of local approximations:
\begin{equation}
u_h = \bigoplus_{k=1}^K u^k_h,  \label{eq:local_approx}
\end{equation}

\noindent
where $u^k_h \in V^k_h$.
Due to \eqref{eq:local_approx}, we can restrict our attention to a single mesh element, dropping the superscript $k$  for simplicity when necessary.
We define the local residual as
\begin{equation*}
R^k_h = \pddt{u^k_h} + \nabla \cdot \vf (u^k_h) - s(u^k_h),
\end{equation*}
and impose that it vanishes locally in a Galerkin sense, i.e.,
\begin{equation*}
\int_{D^k} R^k_h \phi^k_h = 0
\end{equation*}
for any suitably defined test function $\phi^k_h \in V^k_h$.  After  integration by parts, the weak DG formulation is given by:

\begin{equation}
  \int_{D^k} \pddt u^k_h \phi^k_h + 
  \int_{\partial D^k} \vf^* (u_h) \cdot \vn_k \phi^k_h -
  \int_{D^k} \vf (u^k_h) \cdot \nabla \phi^k_h =
  \int_{D^k} s(u^k_h) \phi^k_h    .          \label{eq:dg_parts}
\end{equation}
In Equation \ref{eq:dg_parts}, the physical flux at the element boundary is replaced by a numerical approximation, denoted by $\vf^*$. This guarantees that the flux is single-valued at each edge, enforcing conservation across any edge. Note that all terms in Equation \ref{eq:dg_parts} are local to the $k$-th element, except the numerical flux, which depends on the neighboring elements.  The choice of $\vf^*$ plays a crucial role in the numerical solver's consistency, accuracy and stability.  A popular choice of $\vf^*$ is the Rusanov flux, defined as:
\begin{equation}
\vf^* (u_h) = \vf^* (u^k_h,u^{\bar{k}}_h) =
              \frac{\vf (u^k_h) + \vf (u^{\bar{k}}_h)}{2} -
              \frac{\alpha}{2}(u^{\bar{k}}_h - u^k_h ) \vn_k, \label{eq:rusanov}
\end{equation}
\noindent
where $\bar{k}$ is the index of the neighbor element to $k$ across a given edge, and $\alpha \ge 0$ is a large enough stabilization parameter, usually chosen to be an estimate  of the largest eigenvalue of the hyperbolic system associated to the conservation law.  Finally, $\vn_k$ is the normal unit vector, pointing outwards $D^k$.  
The local solution $u_h^k$ in element $k$ is then written as a linear combination of a polynomial basis $\phi^{(k)}_i$
of $V^k_h:$
    \begin{equation} \label{eq:basis_expansion}
  u_h^k= \sum_{j=1}^{n_{\phi}} \hat u_j^k\phi^{(k)}_j,
    \end{equation}
\noindent
where we have dropped the suffix $h$ in the notation for the polynomial basis. Furthermore, $\hat u_i^k$ denotes the polynomial expansion coefficients, and
$n_{\phi}=(p+1)^2$ represents the cardinality of the polynomial basis, where $p$ represents the maximum degree of the polynomials employed.
Multiple choices exist for the basis set used for the local polynomial spaces. In this study, following e.g., \cite{tumolo:2015}, we employ a modal DG approach, which relies on bivariate Legendre polynomials.
After inserting the basis expansion from Equation \eqref{eq:basis_expansion} in Equation \eqref{eq:dg_parts}
and using the $\phi^{(k)}_i$ as test functions,
 we obtain the following semi-discrete form:
    \begin{equation} \label{eq:ode}
        M \frac{d \mathbf{\hat u}}{dt} =
        \mathbf{h}(\mathbf{\hat u}) \Leftrightarrow \frac{d \mathbf{\hat u}}{dt} = M^{-1}\mathbf{h}(\mathbf{\hat u}).
    \end{equation}
\noindent
Here, the vector $\mathbf{\hat u}$
collects the polynomial expansion coefficients for all elements, and the
matrix $M$ has a block diagonal structure, where the diagonal blocks are the
local mass matrices $M^{(k)}$:
    \begin{equation*} 
        M_{ij}^{(k)} = \int_{D^k} \phi^{(k)}_i \phi^{(k)}_j
    \end{equation*}
associated with each element. Notice that all the terms on the right-hand side have been grouped in the vector function $\mathbf{h}(\mathbf{\hat u}),$
thus obtaining spatial semi-discretization that can be fully discretized by the method of lines approach described below.
Thanks to the use of a DG discretization, the resulting mass matrix can be inverted locally for each element.
Furthermore,  in the case of spherical coordinates, we also simplify the definition of the conserved variables in Equations
\eqref{eq:sph3}
by including the $\cos \theta $ metric terms directly in the  local mass matrix: 
    \begin{equation*}
        M_{ij}^{(k)} = \int_{D^k} \phi_i(\lambda, \theta) \phi_j(\lambda, \theta) \cos(\theta), 
    \end{equation*}
\noindent
so that the conserved variables are given by $h$, $hu$ and $hv$.
 Note that the mass matrices are all identical for a specific longitudinal value.

For the time discretization of Equation~\eqref{eq:ode}, we follow the classical method of lines approach employing Runge-Kutta (RK) methods. More precisely, we use explicit Strong Stability Preserving (SSP) of orders from 1 to 4, denoted later as RK1-RK4, see e.g., \cite{gottlieb:2011}, which can be defined by means of their Butcher tableaux listed in Table \ref{table:Butcher}.
\begin{table}[htb]
    \[
        \begin{array}{c|c}
            0 &  \\
            \hline
            & 1
    
        \end{array}
        \qquad
        \begin{array}{c|c c}
            0 & \\
            1 & 1 & \\
            \hline
            & 1/2 & 1/2
        \end{array}
        \qquad
        \begin{array}{c|ccc}
            0 &  & & \\
            1 & 1 & & \\
            1/2 & 1/4 & 1/4 & \\
            \hline
            & 1/6 & 1/6 & 2/3
        \end{array}
        \qquad
        \begin{array}{c|cccc}
            0 & & & & \\
            1/2 & 1/2 & & & \\
            1/2 & 0 & 1/2 & & \\
            1 & 0 & 0 & 1 & \\
            \hline
            & 1/6 & 1/3 & 1/3 & 1/6
            
        \end{array}
    \]
    \caption{Butcher tableaux of SSP Runge Kutta methods of order 1 to 4}
    \label{table:Butcher}
    \end{table}
These explicit time discretization methods are only conditionally stable. Their stability depends on the value of the non-dimensional parameter known as the Courant number, which is usually defined as 
$ c{\Delta t}/{\cal H}$, where $c$ denotes some estimate of the largest eigenvalue of the underlying hyperbolic system and $\cal H$ is the minimum element diameter. For DG methods and other high-order finite element techniques, however, it is customary
to redefine the Courant number by taking into account the
presence of internal degrees of freedom in each element, see e.g.,
\cite{orlando:2022, orlando:2023, tumolo:2015},
so that a more appropriate definition is in this case
$pc{\Delta t}/{\cal H}$,
where $p$ denotes the maximum element degree. Due to the 
reduction of the effective element size at the poles and the use
of high-order elements, rather small values of the time step
have to be chosen to allow for stable simulations.

\section{Implementations in G4GT and GT4Py} \label{sec:implementation}

We have implemented DG solvers for conservation laws in separate projects with distinct DSLs for planar and spherical geometry.  The Galerkin-for-GridTools (G4GT) and GridTools-for-Python (GT4Py) are both part of the GridTools (GT) \cite{GridTools} ecosystem, which offers an efficient C++ library that is agnostic to the underlying architecture. It makes extensive use of template meta-programming and has backend optimizations for both CPU and GPU architectures. GT was initially designed to target numerical simulations of PDEs using regular grids and finite-difference schemes. Although the latest version of GT supports unstructured grids, the presented implementation relies on its original version.

We remark that G4GT was simply a proof of concept and is {\it no longer supported}. Indeed, the popularity of Python among modern-day programmers led the GT team to switch to the Python-based GT4Py layer, which is actively developed and open source. However, the ideas illustrated in G4GT, in terms of both supported PDE models and computational performance, are educational. Thus, before discussing in detail the GT4Py implementation, which should be regarded as the main tool employed, we emphasize a number of features of the G4GT framework that are complementary to the ones of GT4Py.

\subsection{Galerkin for GridTools (G4GT) implementation}
\label{sec:impl_G4GT}

G4GT is a C++-based extension to the GridTools  library that supports finite element codes.  It relies on GT for the underlying implementation of computation kernels, but also on the Trilinos~\cite{heroux2005overview} libraries Intrepid~\cite{Intrepid} and Epetra~\cite{Epetra}, which provide the numerical support for finite element discretizations and specific linear algebra tools. The G4GT framework provides the link between these libraries, adding a higher-level, user-friendly layer to GT. Additionally, it adds support for finite element discretizations using GT-based codes, which is not present in GT.  

The key steps of the discretization are implemented with GT abstractions, such as the {\tt eval} functor, which subsequently evaluates the element-wise code over the entire domain. Among these steps, the most critical is the computation of the boundary fluxes, as it requires communication between neighboring elements. Its detailed implementation for the horizontal direction is reported below:

\begin{cppcode}
     struct Rusanov_lr{

     using u=gt::accessor<0, enumtype::in, gt::extent<>, 5>;
     using fun=gt::accessor<1, enumtype::in, gt::extent<>, 6>;
     using alpha=gt::global_accessor<2>;
     using normals=gt::accessor<3, enumtype::in, gt::extent<>, 6>;
     using out=gt::accessor<4, enumtype::inout, gt::extent<>, 5>;
     using arg_list=boost::mpl::vector<u, fun, alpha, normals, out>;

     template <typename Evaluation>
     GT_FUNCTION
     static void Do(Evaluation & eval, x_interval) {

     uint_t const num_cub_points=eval.template get_storage_dim<3>(u());
     uint_t const beta_dim=eval.template get_storage_dim<4>(fun());

     for (uint_t face_ : {1,3}) {

         short_t opposite_i = (short_t)(face_==1)?1:(face_==3)?-1:0;
         short_t face_opposite = (short_t)(face_==1)?3:1;
         float_type coeff = -eval(alpha());

         for (short_t qp=0; qp<num_cub_points; qp++) {
             float_type inner_prod1=0.;
             float_type inner_prod2=0.;
             for (uint_t dim=0; dim<beta_dim; dim++){
                 inner_prod1 += eval(fun(0,0,0,qp,dim,face_) * normals(0,0,0,qp,dim,face_)); 
                 inner_prod2 += eval(fun(opposite_i,0,0,qp,dim,face_opposite) * normals(0,0,0,qp,dim,face_));
             }
             eval(out(0,0,0,qp,face_))=(inner_prod1+inner_prod2)/2. - eval(coeff*(u(0,0,0,qp,face_)-u(opposite_i,0,0,qp,face_opposite))/2.);
         }
     }
     }

     };
\end{cppcode}

The main elements of the class can be deduced in a straightforward way.
We simply mention that it relies extensively on the \texttt{gt::accessor} class, and that one should specify whether the accessor is an input or output argument. As the computation is done for each mesh element, for the first three dimensions relative offsets are used as indices.

C++ templating allows for this expansion by making stages of the flux computation on the appropriate grid:

\begin{cppcode}
     auto coords_lr=gt::grid<axis>({1u, 1u, 1u, (uint_t)d1-2u, (uint_t)d1},{0u, 0u, 0u, (uint_t)d2-1u, (uint_t)d2});
     coords_lr.value_list[0] = 0;
     coords_lr.value_list[1] = d3-1;

     auto fluxes_lr=gt::make_computation< BACKEND >(domain_iteration, coords_lr, gt::make_multistage ( execute<forward>(), gt::make_stage< functors::Rusanov_lr > (it::p_u_t_phi_bd(), it::p_fun_bd(),it::p_alpha(), it::p_normals(),  it::p_Rus()) ) );
\end{cppcode}

As the main subject of this work is the GT4Py implementation, we do not delve into additional technicalities of G4GT. However, a more detailed discussion can be found in \cite{discacciati2018implementation}.

\subsection{GT4Py Implementation}  \label{sec:impl_GT4Py}

\begin{figure}[htbp]
        \centering
        \includegraphics[width=0.8\textwidth]{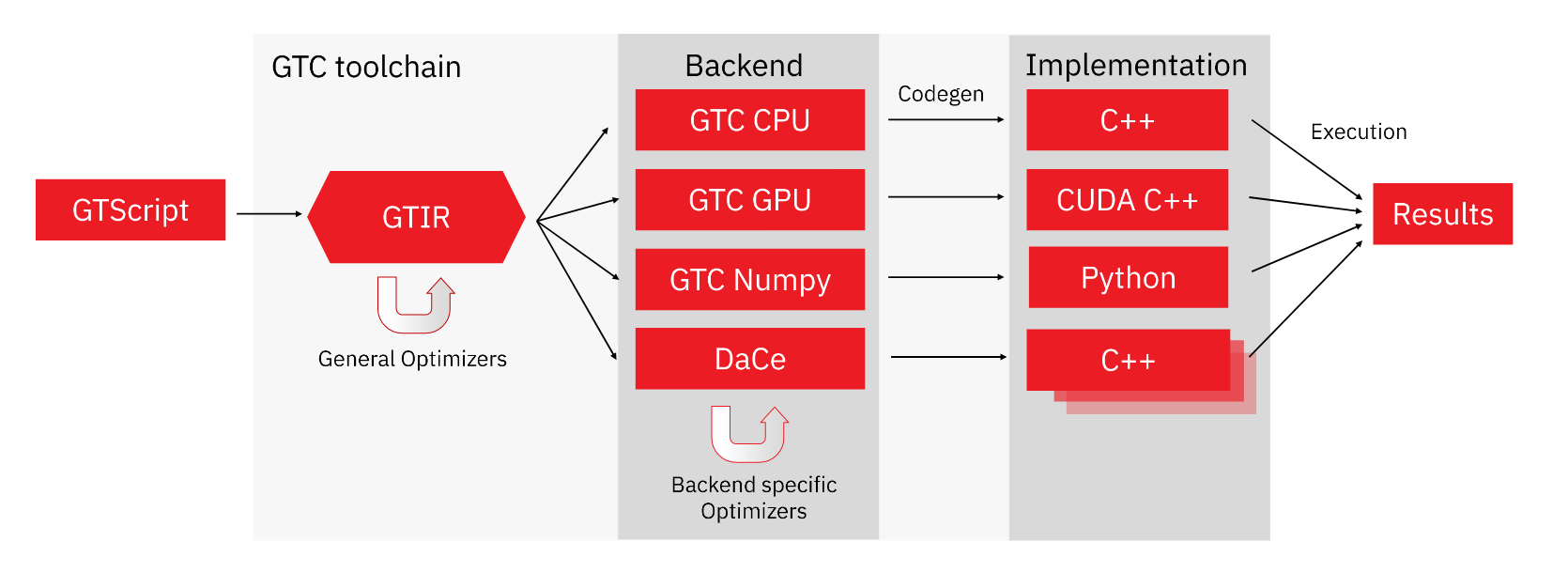}
        \caption[gt4py-pipeline]{GT4Py compilation pipeline. Figure thanks to Till Ehrengruber, CSCS.}
        \label{fig:pipeline}
    \end{figure}

The pipeline of GT4Py is illustrated in Figure \ref{fig:pipeline}. The domain scientist expresses the stencils in a user-friendly Python syntax called GTScript, and this code is then processed through a series of toolchains that applies optimizations and generates a high-performance executable targeting a specific architecture.

\subsubsection{Backends}
GT4Py can compile using various backends; see Table \ref{table:backends} for a complete list of the supported ones at the time of the evaluation. Several others are planned or under development.
    \begin{table}[htb]
        \centering
        \begin{tabular}{cc}\toprule
            Framework & Name \\\midrule
            \multirow{3}{*}{GridTools} & \texttt{gt:cpu\_ifirst} \\
            & \texttt{gt:cpu\_kfirst} \\ 
            & \texttt{gt:gpu} \\
            \midrule
            \multirow{2}{*}{DaCe} & \texttt{dace:cpu} \\
            & \texttt{dace:gpu} \\
            \midrule
            & \texttt{cuda} \\
            & \texttt{numpy} 
            \\\bottomrule
        \end{tabular}
        \caption{List of supported GT4Py backends}
        \label{table:backends}
    \end{table}
Three of the seven backends compatible with GT4Py rely on the GT framework to compile and optimize the stencil computations. They are all characterized with the prefix \texttt{gt:}. The \texttt{gt:cpu\_ifirst} and \texttt{gt:cpu\_kfirst}  both target the CPU architecture, while the \texttt{gt:gpu} backend produces code for the GPU.
In addition, two backends utilize the Data Centric  (\href{https://github.com/spcl/dace}{DaCe}) parallel programming framework \cite{ben-nun_stateful_2019} developed by the Scalable Parallel Computing Lab at Swiss Federal Institute of Technology Zurich, namely \texttt{dace:cpu} and \texttt{dace:gpu} targeting CPUs and GPUs, respectively. At the time of the DG-GT4Py implementation, only prototype implementations of these backends were available, and thus we decided not to include it in the subsequent performance evaluation. 

Alternatively, there is a naive CUDA backend which only utilizes GT utilities, but not its DSL. 
Finally, a NumPy back end exists, which can be used to inspect the generated code for debugging purposes.

\subsubsection{Stencils}
Stencils are special GT4Py functions that operate on fields in a specific domain. Fields store the values of variables at each grid point of the domain.

\paragraph{Declaration}
In the following example, we compute the discretized 2-dimensional Laplacian operator:
    \begin{equation*}
        (\Delta u)_{i,j} = - 4 u_{i,j} + u_{i+1,j} + u_{i-1,j} + u_{i,j+1} + u_{i,j-1} 
    \end{equation*}
which can be written as the following stencil in GT4Py:
   \begin{pythoncode}
    import numpy as np
    import gt4py.gtscript as gtscript
    @gtscript.stencil(backend=backend, **backend_opts)
    def laplacian(
        field: gtscript.Field[np.float64],
        out: gtscript.Field[np.float64]
    ):
        with computation(PARALLEL), interval(...):
            out = - 4 * field + (field[-1,0,0] + field[1,0,0]  + field[0,1,0] + field[0,-1,0])
    \end{pythoncode}
    
In the function decorator, we provide the target backend as well as potential back-end options. The function expects fields as arguments, on which the stencil computations are executed. GT4Py uses the Python-type hinting system to specify the data type of each field, which in this case is \texttt{np.float64}.

The body of the function requires two context managers which define the execution of the stencil in the vertical direction, the first being \texttt{computation} which accepts the arguments \texttt{PARALLEL}, \texttt{FORWARD} or \texttt{BACKWARD}. This defines the scheduling of the execution stencil. The keyword \texttt{PARALLEL}, which we use exclusively for our implementation, indicates that there is no dependence between subsequent vertical levels, and hence they can all be solved in parallel. The keywords \texttt{FORWARD} and \texttt{BACKWARD} define this dependence and indicate the direction in which the vertical levels must be solved. The second context manager \texttt{interval} allows the user to specify the vertical indices for which the stencil will be applied. The `\texttt{...}' is a shorthand notation to select the entire vertical domain.

Finally, we note that the stencil computation is applied for each grid point; hence, relative offsets are used as indices. Note that, if omitted, the offset is assumed to be \texttt{[0,0,0]}.

\paragraph{Invocation}

The above \texttt{laplacian} function can be called using the following command:

    \begin{pythoncode}
    nx, ny, nz = field.shape
    origins = {"field":(1,1,0),"out":(1,1,0)} # or {"_all_":(1,1,0)}
    laplacian(field, out, origin=origins, domain=(nx-2, ny-2, nz))
    \end{pythoncode}

We provide the fields relevant to the stencil computation as arguments to the function. In addition, we add two optional keyword arguments, namely \texttt{domain}, which specifies the domain of execution of the stencil, and \texttt{origin}, which defines the origin for each field. In the case of the \texttt{laplacian} stencil, we set these to ensure that the stencil only operates on the inner part of the domain.

The \texttt{origin} argument indicates relative offset between the different fields. The keyword \texttt{\_all\_} can be utilized to set the same origin for all fields that have not been specified separately. Note that the keyword names inside the \texttt{origins} dictionary refer to the names of the fields in the stencil definition and not to the names of the fields in the call to the stencil. Upon invocation of a stencil, GT4Py searches for a cached version and relies on just-in-time (JIT) compilation in case none is found.

\subsubsection{Storages}

In GT4Py, fields are variables on which stencils can be applied. They store values at each grid point inside a 3-dimensional domain. The DSL provides a storage format for these fields which is a wrapper over the array types  \texttt{{numpy/cupy}.ndarrays} called \texttt{gt4py.storages}, which ensures that the memory layout of the data is compatible with the requested backend. The interface provides several methods for instantiating storages, including \texttt{empty()}, \texttt{ones()} and \texttt{zeros()}, as well as directly from an existing NumPy array using \texttt{from\_array()}. All of these functions require several additional parameters: \texttt{shape} defines the size of the storage in the three dimensions, and \texttt{default\_origin} specifies the default origin to be used in case none is specified during a stencil call. Finally, \texttt{dtype} not only defines the data type of the field but can also be used to assign higher-order tensors to each grid point instead of simple scalar values. These fields are subsequently referred to as higher-dimensional fields.

In the example below, each grid point stores a matrix of size 3x2:
    \begin{pythoncode}
    u = gt4py.storage.zeros(
            backend=backend, default_origin=(1,1,0),
            shape=(4, 4, 2), dtype=(np.float64, (3,2))
        ) 
    \end{pythoncode}
\noindent
Moreover, suppose a field has identical values along one or more spatial dimensions. In that case, GT4Py provides a feature called `masking', which avoids the storage of unnecessary copies of the identical values while still giving the appearance of a full 3-dimensional field.  This can lead to a substantial reduction in memory consumption, which is crucial for large problem sizes.
The previous field can be masked in the vertical direction using:
    \begin{pythoncode}
    u = gt4py.storage.zeros(
            backend=backend, default_origin=(1,1),
            shape=(4, 4), dtype=(np.float64, (3,2)),
            mask=[True, True, False]
        ) 
    \end{pythoncode}

Note that when using a GPU backend, the fields need to be explicitly synchronized from the device back to the host to obtain the results of a stencil computation. In addition, it is recommended to cast the \texttt{gt.storage} to a \texttt{numpy.ndarray} to ensure that the data has indeed been copied from the device. This can be accomplished with the following code snippet:
    \begin{pythoncode}
        x_gt.device_to_host()
        x_np = np.asarray(x_gt)
    \end{pythoncode}

\subsubsection{Frontend}

In this section, we describe the structure of the GT4Py frontend and our contribution to expanding the functionality of higher-dimensional fields.

\paragraph{Abstract Syntax Tree (AST)}

The Python language uses an interpreter which converts the source code of a program into a representation called an Abstract Syntax Tree (AST) before compiling the program to bytecode which is executed by the computer. As the name suggests, the AST represents the logic of the program as a tree structure stripped of the specific syntax used in the source code. This representation provides a versatile way to inspect and modify Python applications.

In the case of GT4Py, the frontend parses the Python AST and converts it into a series of custom ASTs through the pipeline (Figure \ref{fig:pipeline}), which provide additional information necessary for the backends to produce well-optimized executables.

\paragraph{Limited support for higher-dimensional fields}

For our implementation, we represent each DG element by a grid point in GT4Py. Each grid point is thus assigned a vector that stores its polynomial expansion coefficients and hence yields a higher-dimensional field as a data structure. We refer to this additional dimension of the field as \texttt{data\_dims}.

Initially, the support for these vector-valued fields in GT4Py was limited.
In particular, there was no functionality for performing element-wise operations between fields with respect to the \texttt{data\_dims} dimension.
Indeed, these operations needed to be explicitly written out for each vector component, reducing their utility to scalar fields.
The following example illustrates a stencil performing an element-wise multiplication between two higher-dimensional fields:
    \begin{pythoncode}
    @gtscript.stencil(backend=backend)
    def mult(
        field1: gtscript.Field[(np.float64, (3,))],
        field2: gtscript.Field[(np.float64, (3,))],
        out: gtscript.Field[(np.float64, (3,))]

    ):
        with computation(PARALLEL), interval(...):
            out[0,0,0][0] = field1[0,0,0][0] * field2[0,0,0][0]
            out[0,0,0][1] = field1[0,0,0][1] * field2[0,0,0][1]
            out[0,0,0][2] = field1[0,0,0][2] * field2[0,0,0][2]
    \end{pythoncode}
The first set of indices represents the relative offsets between the fields, while the second set of indices refers to the actual components of the \texttt{data\_dims} dimension. Note that in this case, the relative offsets cannot be omitted and need to be specified explicitly.

\paragraph{Loop unroller}

We have implemented automatic element-wise operations for vector-valued fields in order to facilitate their use.
This has been accomplished by modifying an intermediate GT4Py AST called \texttt{definition\_ir}.
In this representation, we have access to the size of \texttt{data\_dims} for each field which is essential to determine if an operation between two fields should be performed element-wise.
Our contribution allows rewriting the previous multiplication stencil using the following simple syntax:
    \begin{pythoncode}
    # ...
    with computation(PARALLEL), interval(...):
        out = field1 * field2
    \end{pythoncode}

The goal is to modify the GT4Py frontend so that the code snippet above produces the same AST as the previous explicitly unrolled stencil. To implement this functionality, we created two helper classes which apply the necessary transformations to the \texttt{definition\_ir}. 
The first one, called \texttt{UnRoller}, receives as argument an AST node representing the right-hand side  of an assignment operation and returns a list of AST nodes where each element of the list pertains to a different index of the higher-dimensional field.
The second one is called \texttt{UnVectorisation}. It invokes the \texttt{UnRoller} and checks that the dimensions of the returned list match the dimensions of the field on the left-hand side of the assignment operation. If so, it creates the list of AST assignment nodes with the corresponding indices.

Our implementation supports not only chaining together multiple operations on higher-dimensional fields but also broadcasting of scalar values for scalar-vector operations. The syntax and functionality should be intuitive for anyone familiar with the NumPy package.

\paragraph{Matrix multiplication}

An additional operation that we required for our DG scheme was a matrix-vector multiplication between higher-dimensional fields. This was incorporated into our existing framework and can be invoked using the "\texttt{@}" operator. Also, the multiplication of a vector by the transposed of a matrix can be achieved by appending the matrix with the "\texttt{T}" attribute. This leads to the following syntax:
\begin{pythoncode}
    @gtscript.stencil(backend=backend)
    def matmul(
        matrix: gtscript.Field[(np.float64, (3, 2))],
        vec: gtscript.Field[(np.float64, (3,))],
        out: gtscript.Field[(np.float64, (2,))]
    )
        with computation(PARALLEL), interval(...):
            out = matrix.T @ vec
\end{pythoncode}

\paragraph{DG solver: precomputation}

At the start of the execution of the program, the GT4Py solver precomputes on the CPU certain variables  that remain constant during the whole simulation. This includes the computation of the inverse mass matrix, as well as the Gauss-Legendre quadrature points and weights for numerical integration. A helper class called Vander, defined in \texttt{vander.py}, contains all the Vandermonde matrices  required to evaluate the polynomials stored as modal expansion coefficients at nodal values in the domain (see e.g., \cite{hesthavenNodalDiscontinuousGalerkin2008}). These matrices are instantiated as fields using \texttt{gt4py.storages}.

\paragraph{DG solver: stencils}

All subsequent computations are carried out using stencils in GT4Py. An example stencil is presented subsequently, related to our DG solver. Applying the theory derived in Section \ref{sec:math} for the linear, constant-coefficient advection problem

\begin{equation} \label{eq:lin_adv}
        \frac{\partial u}{\partial t} + \nabla \cdot (\boldsymbol{\beta} u) = 0,
    \end{equation}   
    with e.g., $\boldsymbol{\beta} = [1, 1]^T,$
we will need to evaluate an integral of the following form:
    \begin{equation*}
        \int_{D^k} \left [\beta_1 u \frac{\partial \phi}{\partial x} + \beta_2 u \frac{\partial \phi}{\partial y} \right ] dx dy.
    \end{equation*}
This integral can be computed using the stencil below:
\begin{pythoncode}
    #...
    with computation(PARALLEL), interval(...):
        u_qp = phi @ u_modal
        fx = u_qp * 1
        fy = u_qp * 1
        rhs = determ * (phi_grad_x.T @ (fx * w) / bd_det_x
                        + phi_grad_y.T @ (fy * w) / bd_det_y)
\end{pythoncode}
In line 3, the modal expansion coefficients are mapped to nodal values at the quadrature points. In lines 4 and 5, the flux function in the x and y directions is applied. In this simple case, the flux function is the identity due to the constant velocity field $\boldsymbol{\beta} = [1,1]^T$. Finally, in lines 6 and 7, the numerical integration is performed. The scalar field \texttt{w} represents the quadrature weights while the matrix-valued field \texttt{phi\_grad\_x/y} contains the spatial derivatives of the basis functions. The terms \texttt{determ}, \texttt{bd\_det\_x/y} denote the Jacobians arising from the mapping of the element in physical space onto a reference element. Although not reported here, this description easily generalizes to the SWE case.

\section{Code validation}  \label{sec:validation}

Firstly, several unit tests have been performed to assess  the correctness of our implementations. For the GT4Py implementation, they rely on the existing testing infrastructure of GT4Py, which verifies the success of the code generation as well as the code execution on all backends when compared with a reference Numpy implementation. 

Subsequently, both the G4GT and GT4Py implementations were validated on benchmarks derived from the shallow water test suite~\cite{williamsonStandardTestSet1992}, as well as on tests on a planar geometry presented in \cite{tumolo:2013}. These include a convergence test on linear advection of a smooth profile and on a geostrophic zonal flow, the simulation of geostrophic adjustment on the plane, and that of a Rossby-Haurwitz wave in spherical geometry.
Although our final goal is the simulation of the SWE on the sphere, all the presented tests provide useful insight from a numerical point of view and contribute to a progressive increase in the complexity of the solutions.

\subsection{Linear advection convergence of a smooth initial condition}

Since both implementations are essentially solving the same problem, albeit with slightly different flux calculations and numerical implementations, we summarize the convergence results for both in this section.
Specifically, we  apply   our  DG discretization to  planar linear advection on the unit square, assuming periodic boundary conditions and a constant velocity field $\boldsymbol{\beta}, $
see Equation \eqref{eq:lin_adv}.
We consider a smooth initial condition:
      $   u_0(x, y) = \sin(2 \pi x) \sin(2 \pi y), $
       which allows us to achieve optimal convergence rates. 
The analytic solution of Equation (\ref{eq:lin_adv})  evolves without changing shape in the direction of the velocity field. Due to the periodic boundary conditions, the solution will coincide with the initial condition after one full rotation, i.e., at time $T=1.$
We use a uniform mesh with $K$ elements obtained from the tensor product of $\sqrt{K}$ elements in each of the coordinate directions.
We measure the error of the numerical approximation using the $L^2$ norm at the final simulation time, and we denote it by $\epsilon$.
The expected spatial convergence order for DG methods is given by \cite{hesthavenNodalDiscontinuousGalerkin2008}:
    \begin{equation}
    \label{eq:conv_order}
        \epsilon \sim O(h^{p+1}),
    \end{equation}
\noindent
where $h$ is the characteristic mesh size and $p$ is the degree of the local polynomials.
To estimate the convergence rate, we compute the discretization error using two different meshes with characteristic sizes $h_1, \ h_2,$
that we denote as $\epsilon_1, \epsilon_2, $ respectively.
Then, the estimated rate, denoted by $r$, is computed as
\begin{equation}
r = \frac{\log(\epsilon_1)-\log(\epsilon_2)}{\log(h_1)-\log(h_2)}.
\end{equation}
The results obtained with the G4GT implementation are reported in Table \ref{tab:linear_advection_conv} and agree with the theoretical expectations, thus validating the implementation. Not surprisingly, the GT4Py implementation  achieves nearly identical convergence results to the G4GT version, as shown in   Table \ref{table:sine_conv}.

\begin{table}[htbp]
	\centering
	\begin{tabular}{ccccccc}
		\toprule	
		\multirow{2}{*}{$K$} & \multicolumn{2}{c}{$p=1$} & \multicolumn{2}{c}{$p=2$} &
		\multicolumn{2}{c}{$p=3$} \\ \cmidrule(r){2-3} \cmidrule(lr){4-5} \cmidrule(l){6-7}
		& $\epsilon$ & $r$ & $\epsilon$ & $r$ & $\epsilon$ 
  & $r$ \\ \midrule
		$10^2$ & $\expnumber{1.343}{-2}$ & - & $\expnumber{1.050}{-3}$ & - & $\expnumber{3.780}{-5}$ & - \\ 
		$20^2$ & $\expnumber{3.369}{-3}$ & \textbf{2.00} & $\expnumber{1.329}{-4}$ & \textbf{2.98} & $\expnumber{2.030}{-6}$ & \textbf{4.22} \\
		$40^2$ & $\expnumber{8.405}{-4}$ & \textbf{2.00} & $\expnumber{1.666}{-5}$ & \textbf{3.00} & $\expnumber{1.302}{-7}$ & \textbf{3.96} \\ 
		$80^2$ & $\expnumber{2.099}{-4}$ & \textbf{2.00} & $\expnumber{2.084}{-6}$ & \textbf{3.00} & $\expnumber{8.345}{-9}$ & \textbf{3.96} \\
		$160^2$ & $\expnumber{5.246}{-5}$ & \textbf{2.00} & $\expnumber{2.611}{-7}$ & \textbf{3.00} & & \\ 		
		\bottomrule
	\end{tabular}
	\caption{$L^2$ errors $\epsilon$ and estimated rate of convergence $r$ for the linear advection problem, G4GT implementation.}
	\label{tab:linear_advection_conv}
\end{table}

   \begin{table}[htb]
        \centering
        \begin{tabular}{ccccccc} 
            \toprule
		\multirow{2}{*}{$K$} & \multicolumn{2}{c}{$p=1$} & \multicolumn{2}{c}{$p=2$} &
		\multicolumn{2}{c}{$p=3$} \\ \cmidrule(r){2-3} \cmidrule(lr){4-5} \cmidrule(l){6-7}
		& $\epsilon$ & $r$ & $\epsilon$ & $r$ & $\epsilon$ 
  & $r$ \\ \midrule
            $20^2$  & 4.204e-3 & - &  1.330e-4 & - & 2.061e-6 & - \\
            $40^2$ & 9.004e-4 & \textbf{2.22} & 1.666e-5 & \textbf{2.99} & 1.288e-7 & \textbf{4.00} \\
            $80^2$ & 2.139e-4 & \textbf{2.07} & 2.084e-6 & \textbf{2.99} & 8.049e-9 & \textbf{4.00} \\
            $160^2$ & 5.212e-5 & \textbf{2.02} & 2.606e-7 & \textbf{3.00} & 5.030e-10 & \textbf{4.00} \\\bottomrule
        \end{tabular}
        \caption{$L^2$ errors $\epsilon$ and estimated rate of convergence $r$ for the linear advection problem, GT4Py implementation.}
        \label{table:sine_conv}
    \end{table}

\subsection{G4GT implementation:  geostrophic adjustment for planar SWE}

For the G4GT implementation, we consider the SWE discretized on a Cartesian mesh in a planar domain. Specifically, we want to show that the scheme is able to reproduce the geostrophic adjustment process, see, for example, the discussion in
\cite{tumolo:2013}. Starting from a perturbation of the equilibrium state corresponding to a constant water height, gravitational and rotational forces interact, so that only part of the energy is transported away from the center, leading to a nontrivial stationary solution profile.
\begin{figure}[htbp]
	\centering
	\begin{subfigure}{0.49\textwidth}
		\centering
		\includegraphics[width=\linewidth]{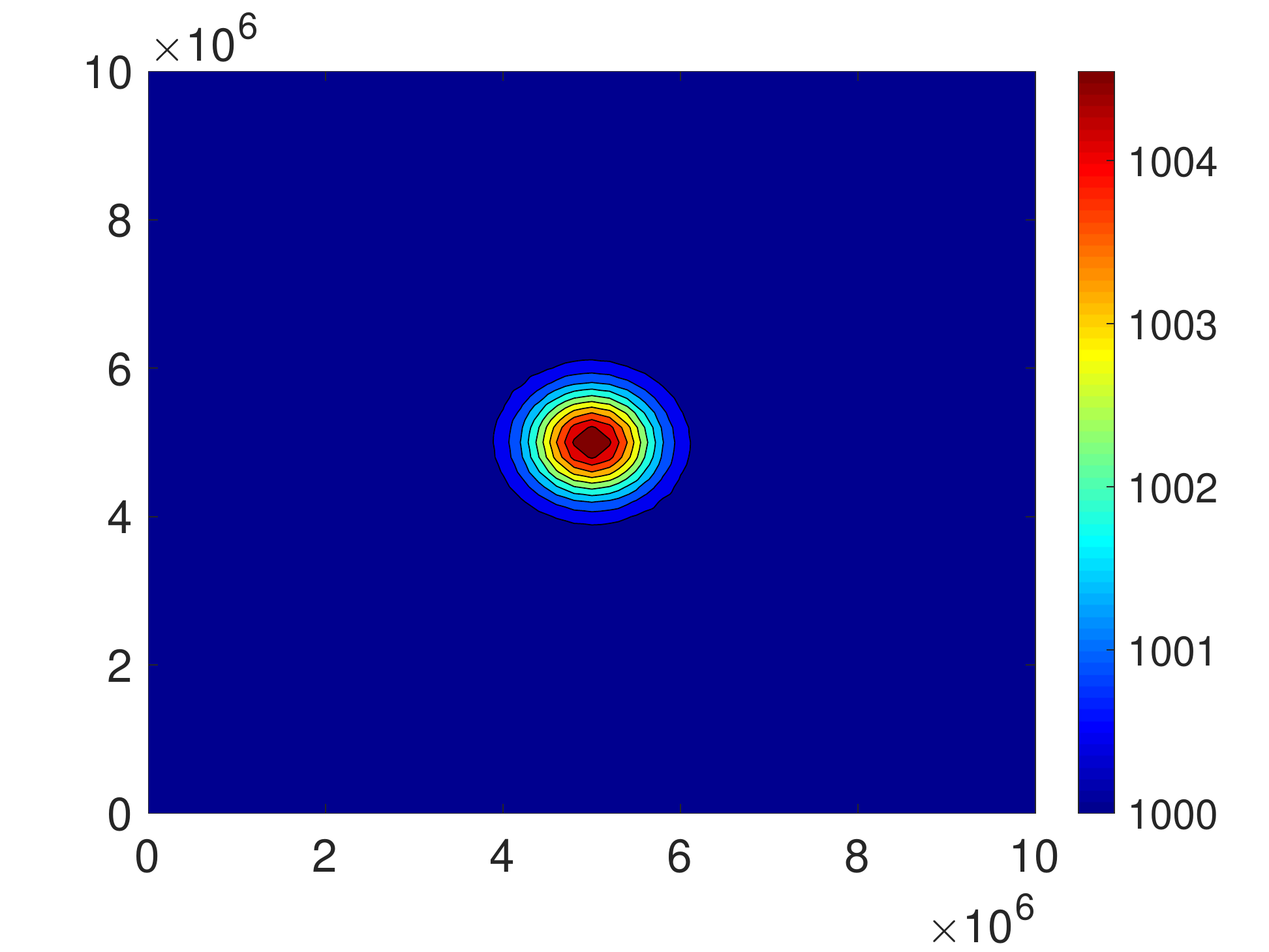}
		\subcaption{Initial height.}
	\end{subfigure}
	\begin{subfigure}{0.49\textwidth}
		\centering
		\includegraphics[width=\linewidth]{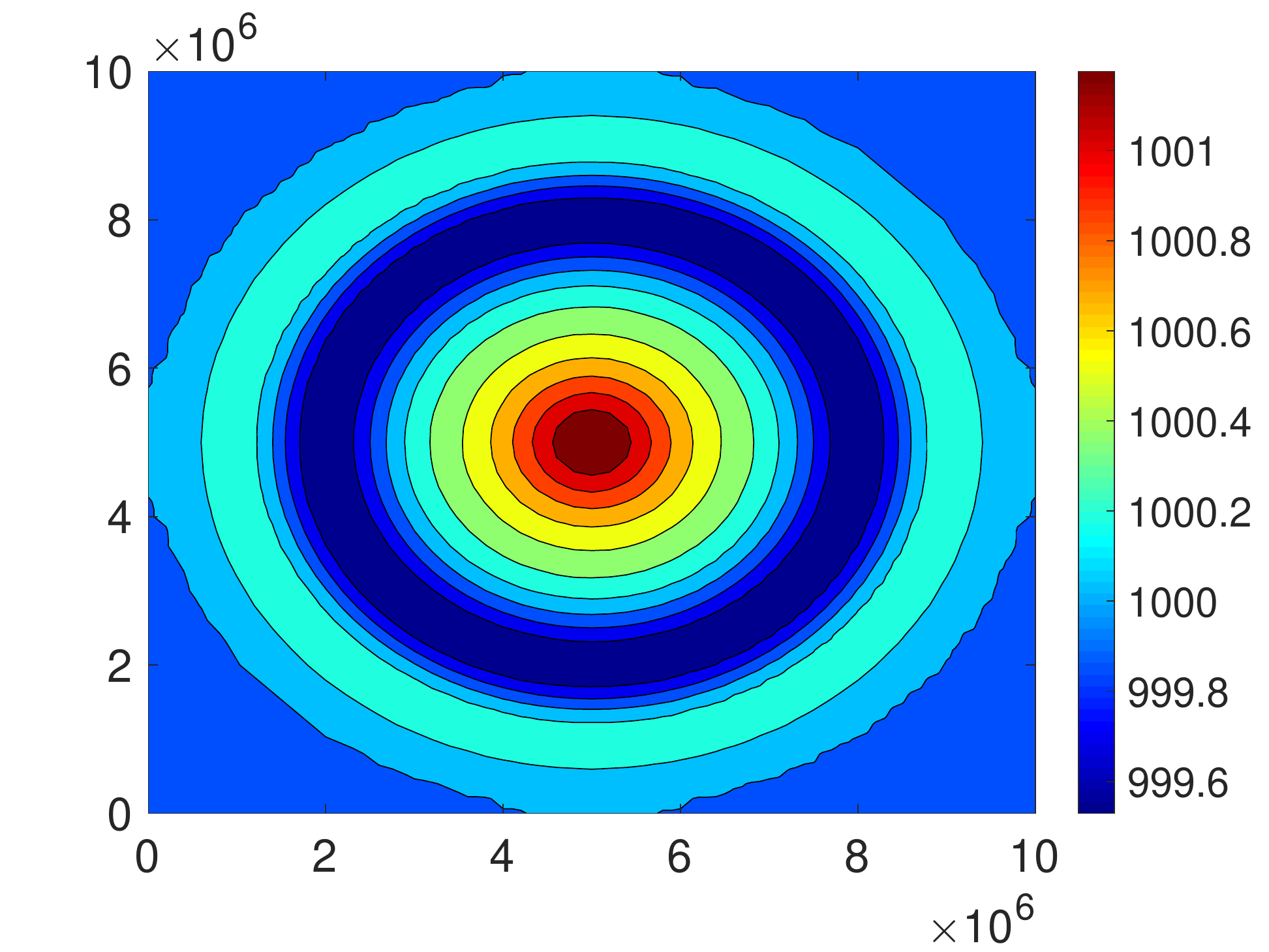}
		\subcaption{Height.}
	\end{subfigure}
	
	\begin{subfigure}{0.49\textwidth}
		\centering
		\includegraphics[width=\linewidth]{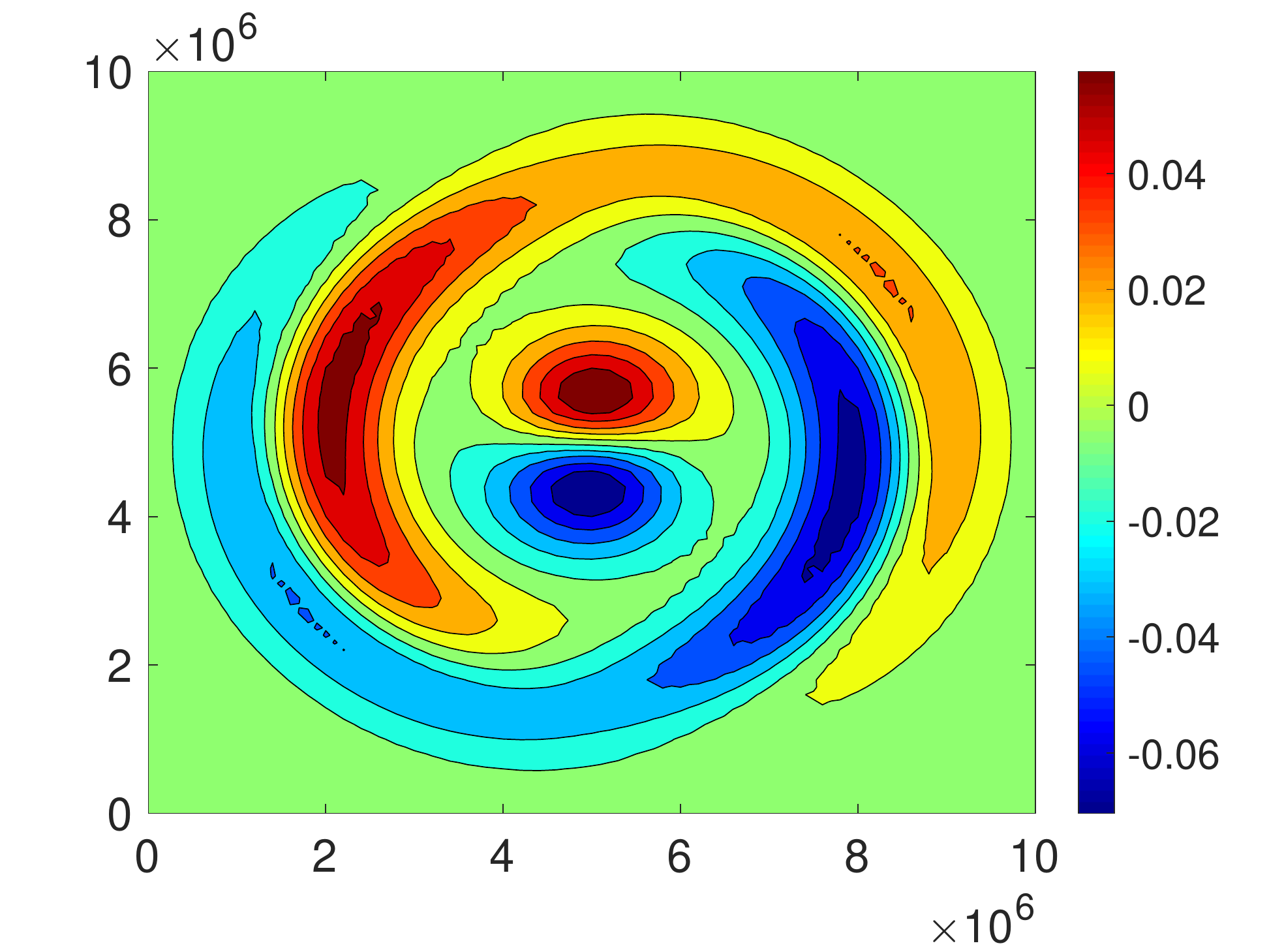}
		\subcaption{Horizontal velocity.}
	\end{subfigure}
	\begin{subfigure}{0.49\textwidth}
		\centering
		\includegraphics[width=\linewidth]{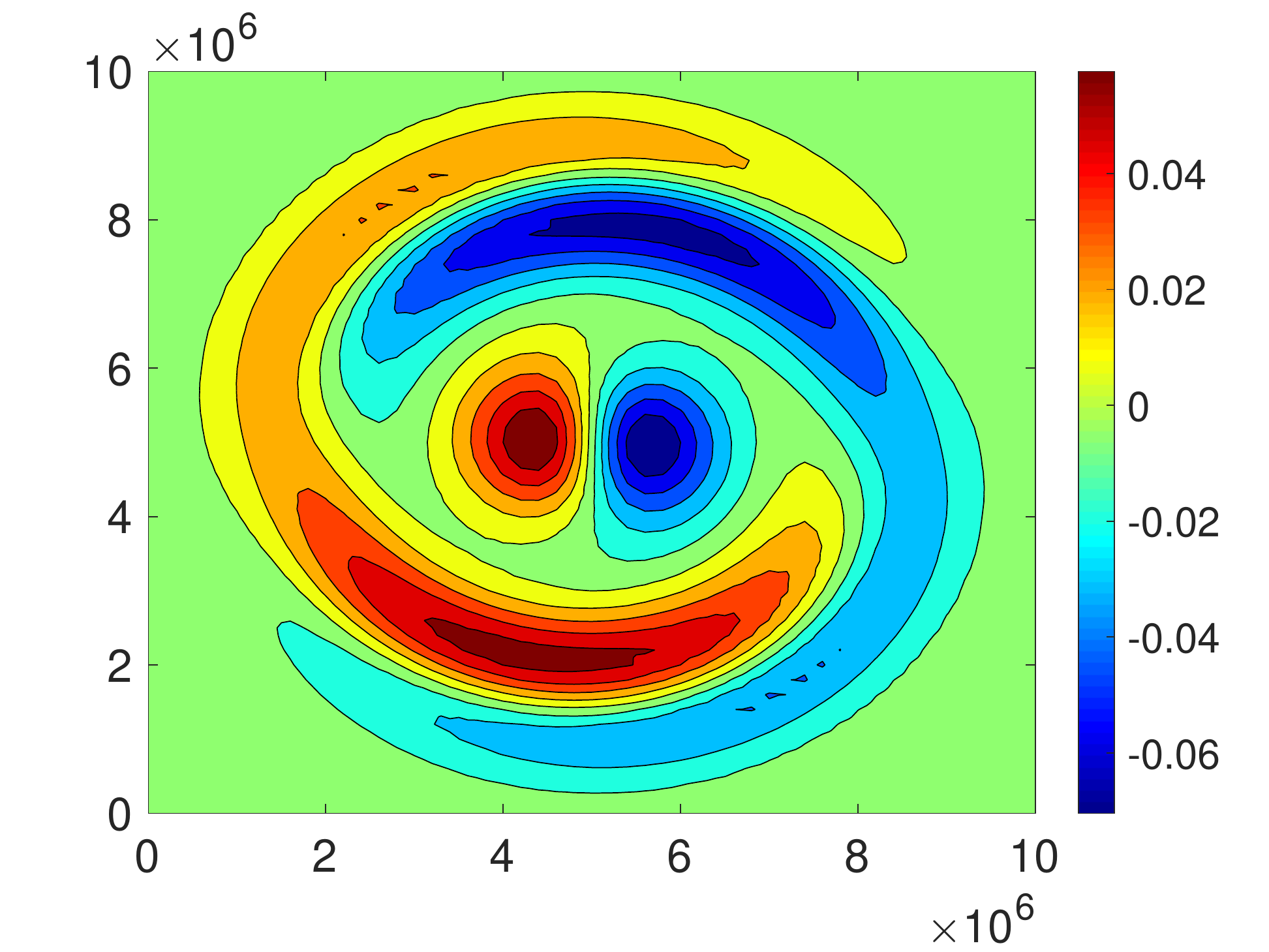}
		\subcaption{Vertical velocity.}
	\end{subfigure}
	
	\caption{Numerical results for the geostrophic adjustment test case.}
	\label{fig:geostrophic_adj}
\end{figure}
\noindent 
Consider a square domain $\Omega=[0,L]^2$ with $L=10^7 \ \rm m$ and a final time of $T=36000 \ \rm s$. The initial velocities and momenta are set to zero, while the height $h$ is equal to
\begin{equation}
    h=h_0+h_1 \exp \left( - \frac{(x-L/2)^2+(y-L/2)^2}{2 \sigma^2}  \right),
\end{equation}
where $h_0=1000 \ \rm m$, $h_1=5 \ \rm m$ and $\sigma=L/20 \ \rm m$. Assuming an $f$-plane approximation, the Coriolis parameter $f$ is chosen to be constant and equal to $10^{-4} \ \rm s^{-1}$. The problem is completed with periodic boundary conditions. The simulation has been run using $50$x$50$ spatial elements, a polynomial degree $r=3$ and the RK4 scheme in time with step $\Delta t=100 \ \rm s$. The results are reported in \figref{fig:geostrophic_adj}.
The solution is consistent with the results reported in \cite{tumolo:2013}. 

\subsection{GT4Py implementation: geostrophic zonal flow for SWE on the sphere}
After validating the planar version of the GT4Py implementation, we consider two of the classical test cases in spherical geometry introduced in \cite{williamsonStandardTestSet1992}
for the shallow water equations.
Periodic boundary conditions were applied in the longitudinal direction, while in the latitudinal direction the fluxes were set to zero. Indeed, since the edges become singular at the poles, the flux through them must be zero.

In the benchmark denoted as test case 2 in \cite{williamsonStandardTestSet1992}, a stationary
zonal flow in geostrophic equilibrium is considered. We perform a convergence test
for the spatial discretization, using for all polynomial degrees the RK4 scheme for time discretization
with  time steps chosen for each resolution so as to keep the Courant number constant and at a very small value. The test case has been run until time $T=2 $ days on meshes of increasing resolutions. The results are reported in Table \ref{tab:test2_conv} and show a 
convergence behavior entirely analogous to that of the linear advection case in planar geometry.

\begin{table}[htbp]
	\centering
	\begin{tabular}{ccccccc}
		\toprule	
		\multirow{2}{*}{$K$} & \multicolumn{2}{c}{$p=1$} & \multicolumn{2}{c}{$p=2$} &
		\multicolumn{2}{c}{$p=3$} \\ \cmidrule(r){2-3} \cmidrule(lr){4-5} \cmidrule(l){6-7}
		& $\epsilon$ & $r$ & $\epsilon$ & $r$ & $\epsilon$ 
		& $r$ \\ \midrule
		$10^2$ & $\expnumber{9.366}{-2}$ & - & $\expnumber{1.020}{-3}$ & - & $\expnumber{6.864}{-5}$ & - \\ 
		$20^2$ & $\expnumber{1.984}{-3}$ & \textbf{5.56} & $\expnumber{1.085}{-4}$ & \textbf{3.23} & $\expnumber{3.951}{-6}$ & \textbf{4.08} \\
		$40^2$ & $\expnumber{4.508}{-4}$ & \textbf{2.13} & $\expnumber{1.490}{-5}$ & \textbf{2.86} & $\expnumber{2.362}{-7}$ & \textbf{4.06} \\ 
		$80^2$ & $\expnumber{1.111}{-4}$ & \textbf{2.02} & $\expnumber{1.986}{-6}$ & \textbf{2.90} & $\expnumber{1.471}{-8}$ & \textbf{4.00} \\ 		
		\bottomrule
	\end{tabular}
	\caption{$L^2$ errors $\epsilon$ and estimated rate of convergence $r$ for the Williamson test case 2, GT4Py implementation.}
	\label{tab:test2_conv}
\end{table}

\subsection{GT4Py implementation: Rossby-Haurwitz wave for SWE on the sphere}

\begin{figure}[htbp]
	\centering
	\begin{subfigure}[t]{0.32\textwidth}
		\centering
		\includegraphics[width=\linewidth]{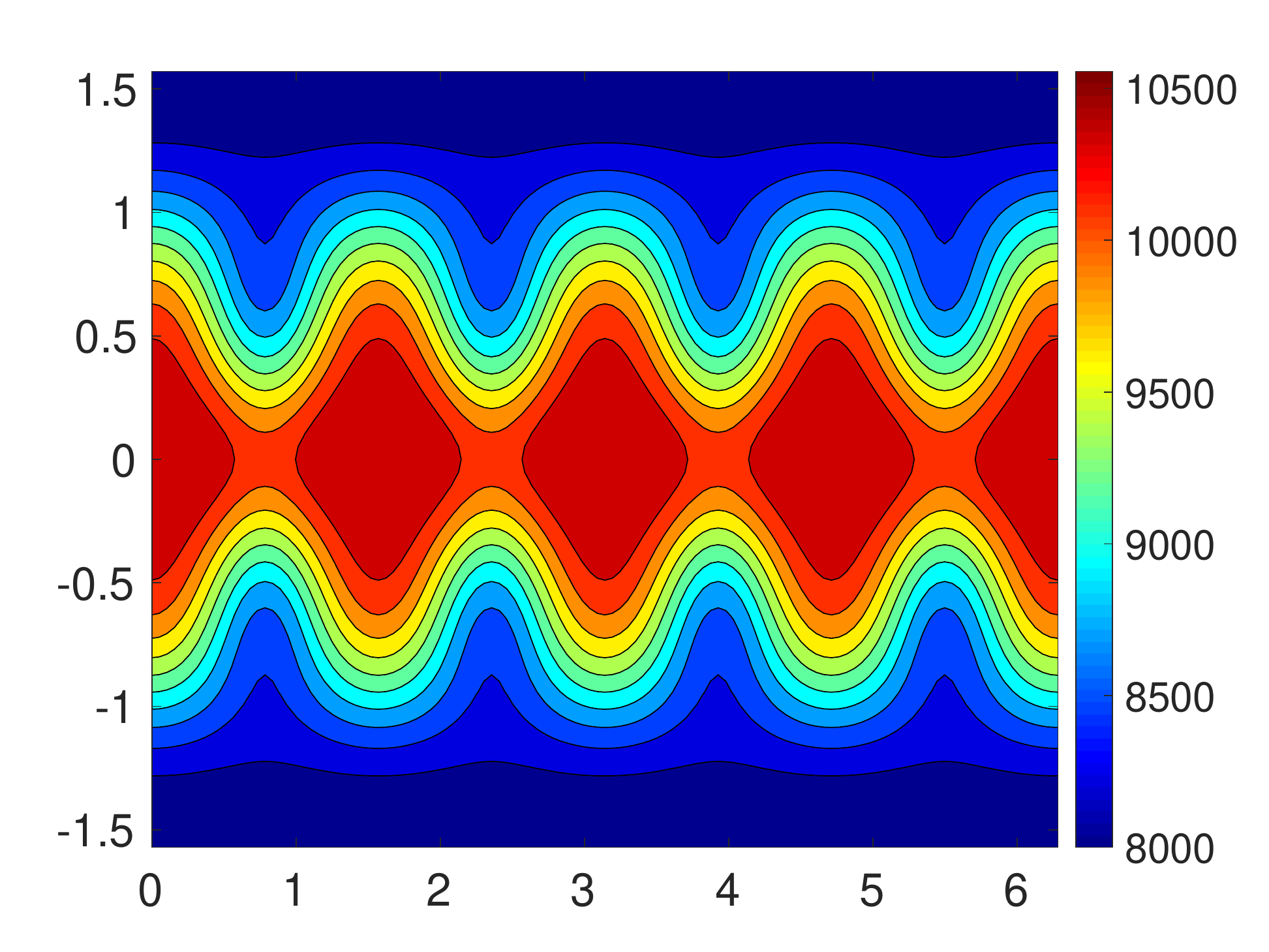}
		\subcaption{$h$ at $t=0$ days.}
	\end{subfigure}
	\begin{subfigure}[t]{0.32\textwidth}
		\centering
		\includegraphics[width=\linewidth]{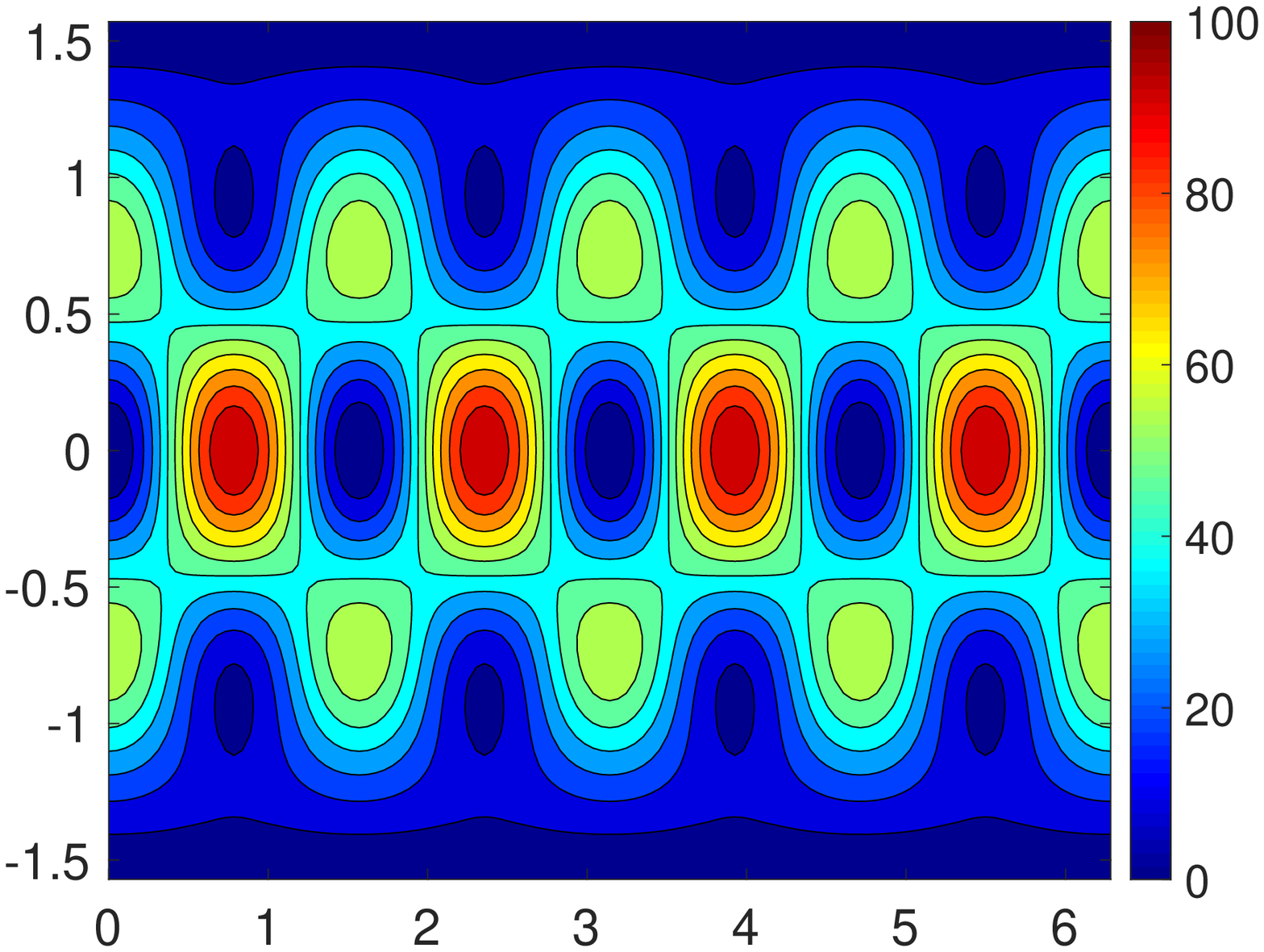}
		\subcaption{$u$ at $t=0$ days.}
	\end{subfigure}
	\begin{subfigure}[t]{0.32\textwidth}
		\centering
		\includegraphics[width=\linewidth]{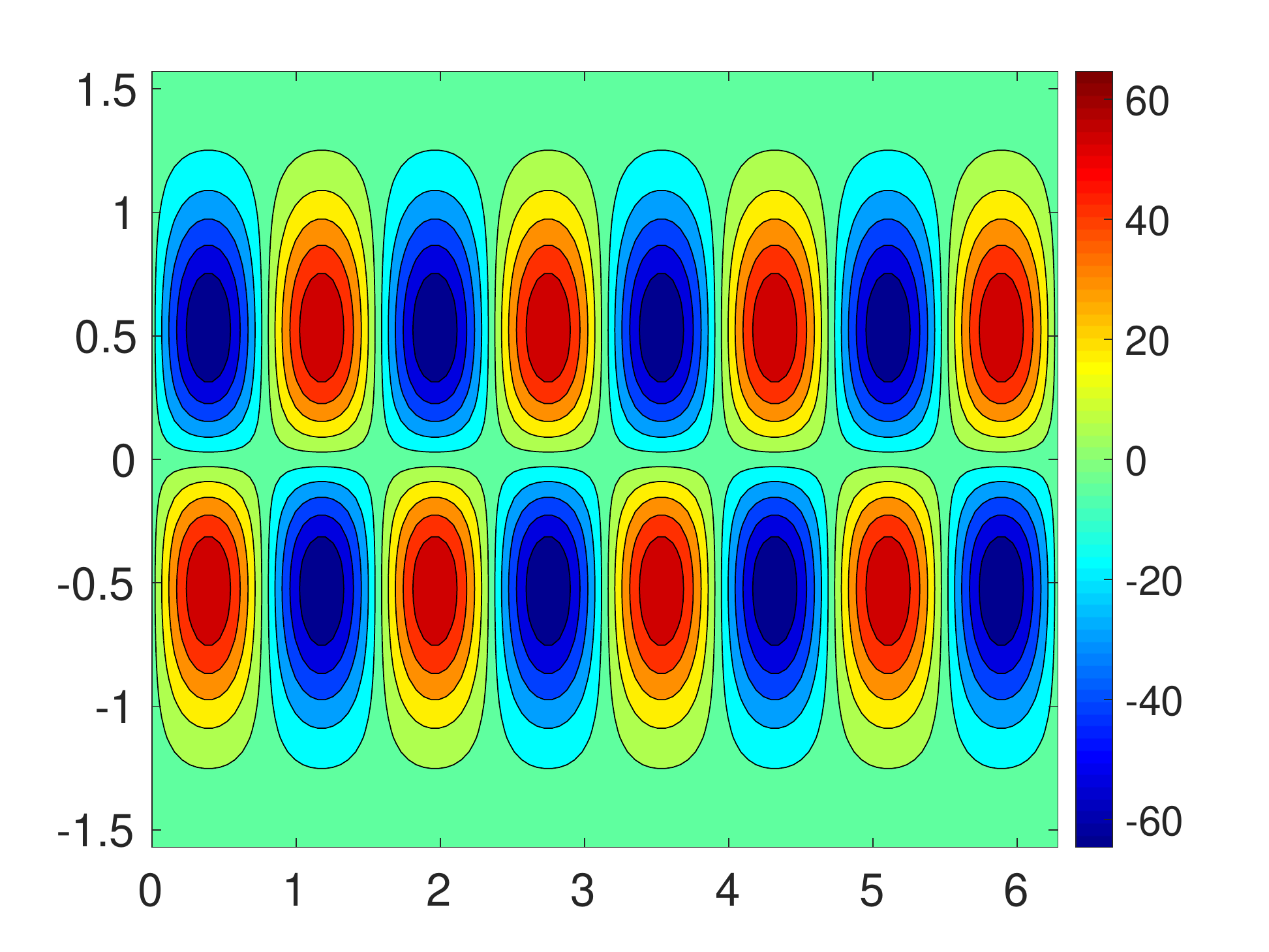}
		\subcaption{$v$ at $t=0$ days.}
	\end{subfigure}

	\begin{subfigure}[t]{0.32\textwidth}
		\centering
		\includegraphics[width=\linewidth]{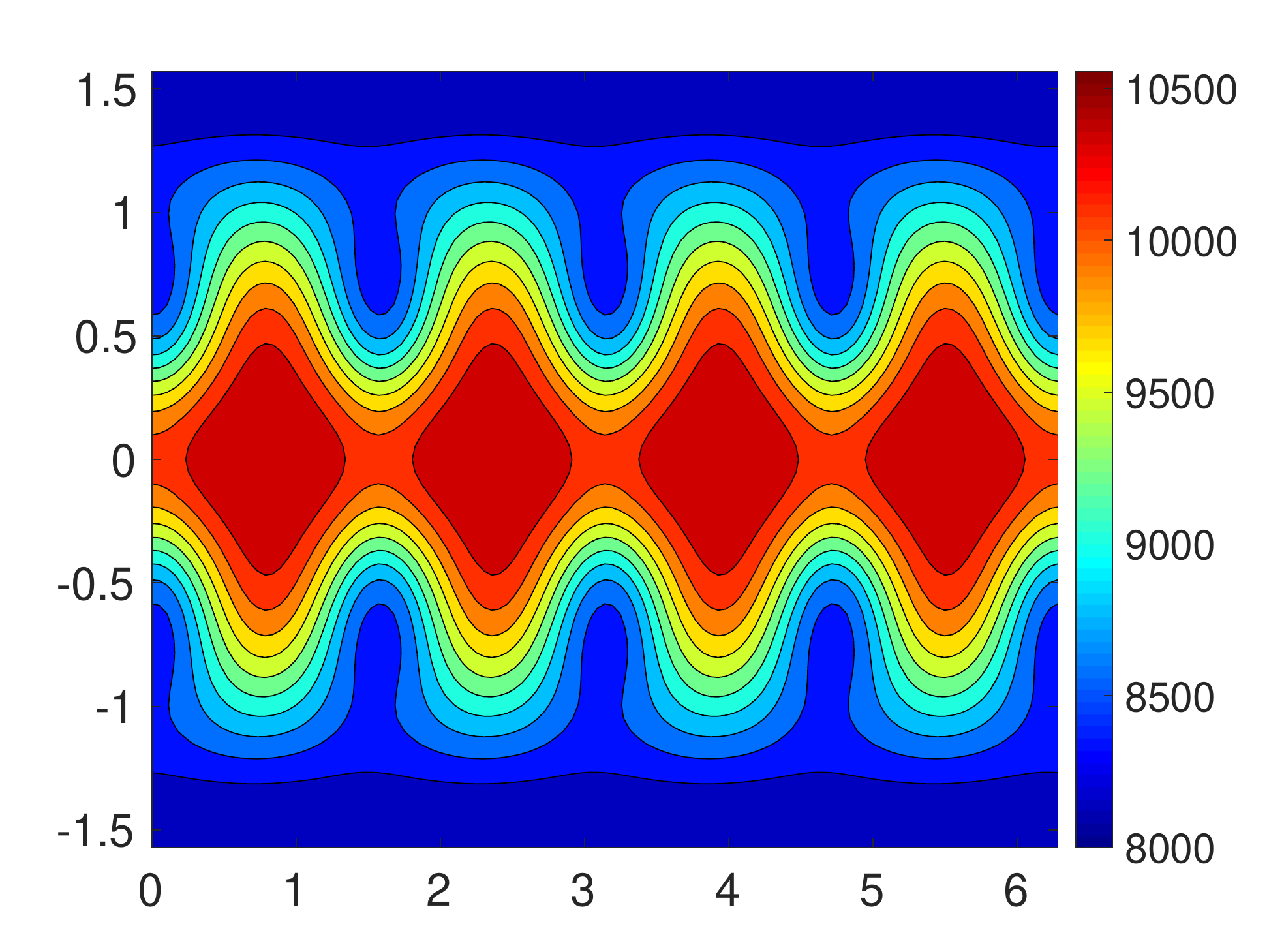}
		\subcaption{$h$ at $t=4$ days.}
	\end{subfigure}
	\begin{subfigure}[t]{0.32\textwidth}
		\centering
		\includegraphics[width=\linewidth]{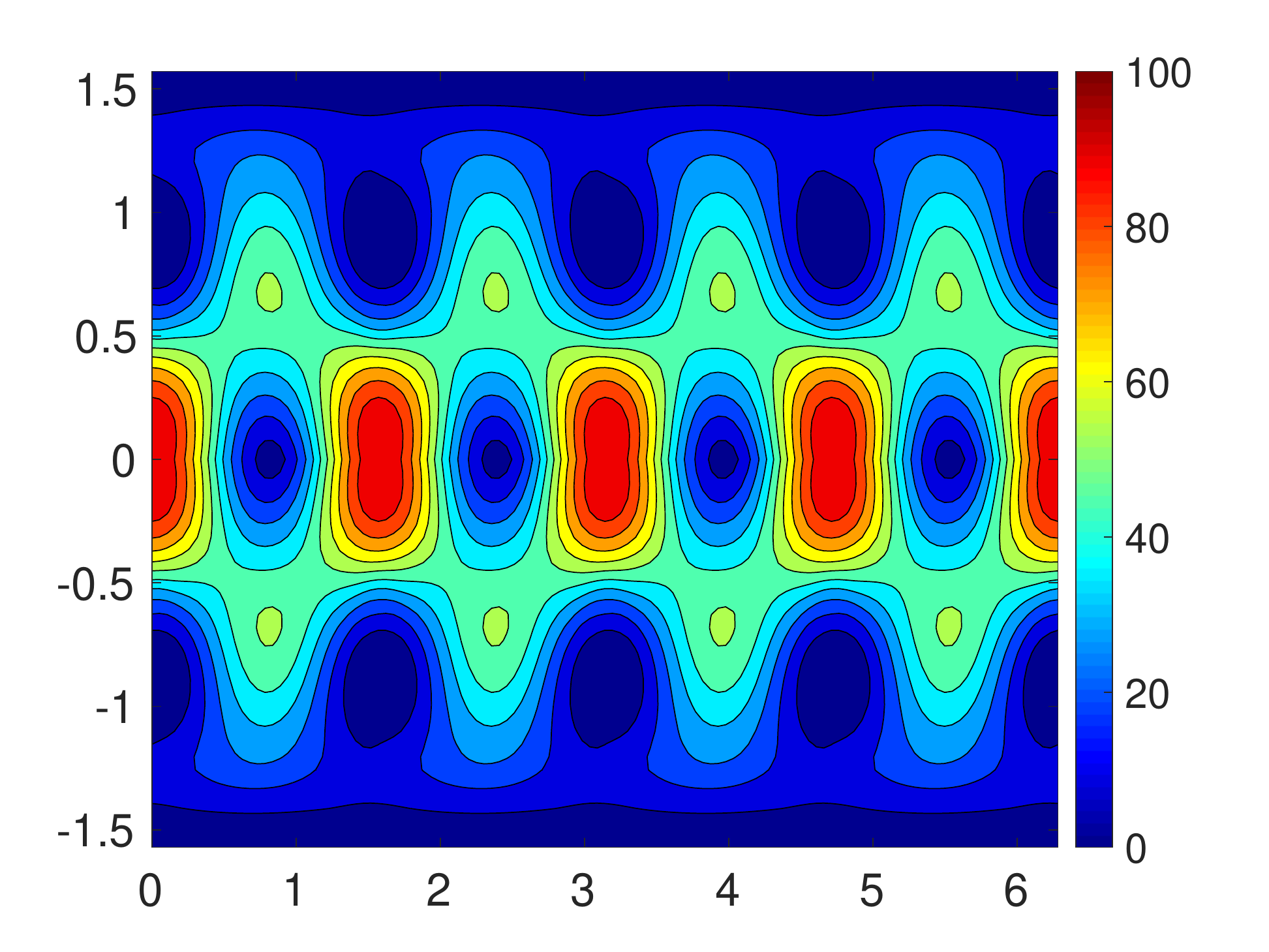}
		\subcaption{$u$ at $t=4$ days.}
	\end{subfigure}
	\begin{subfigure}[t]{0.32\textwidth}
		\centering
		\includegraphics[width=\linewidth]{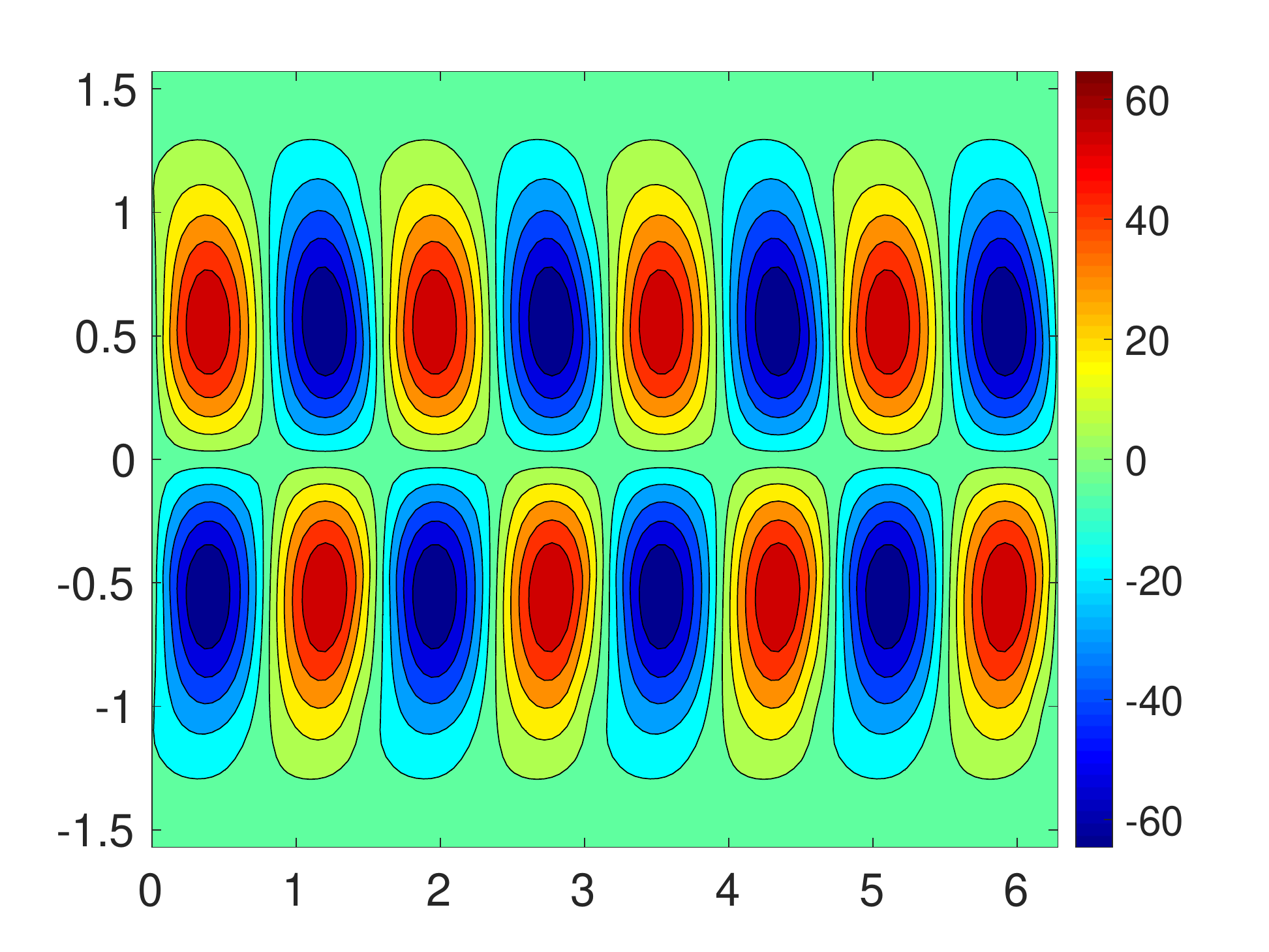}
		\subcaption{$v$ at $t=4$ days.}
	\end{subfigure}

	\begin{subfigure}[t]{0.32\textwidth}
		\centering
		\includegraphics[width=\linewidth]{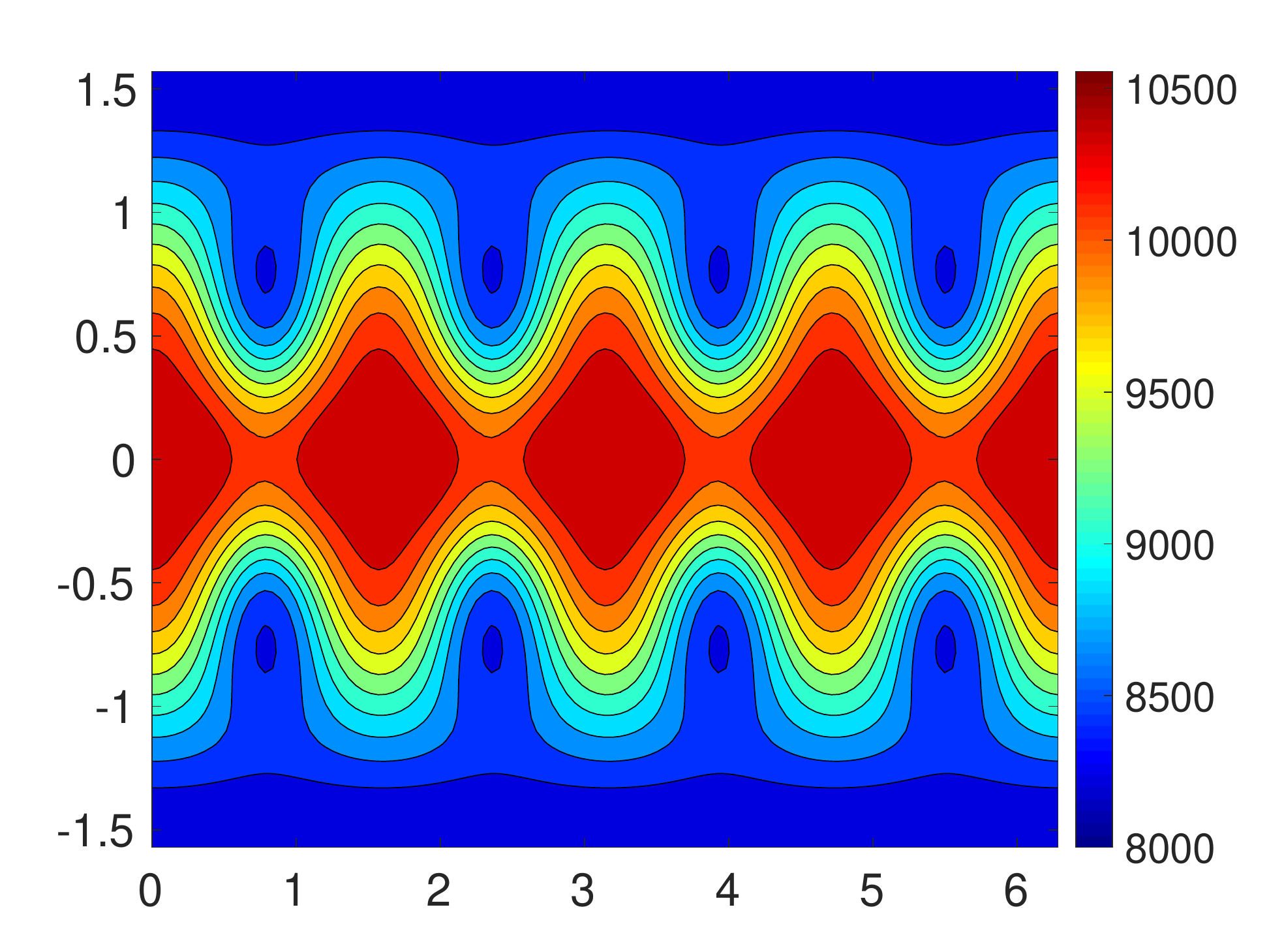}
		\subcaption{$h$ at $t=8$ days.}
	\end{subfigure}
	\begin{subfigure}[t]{0.32\textwidth}
		\centering
		\includegraphics[width=\linewidth]{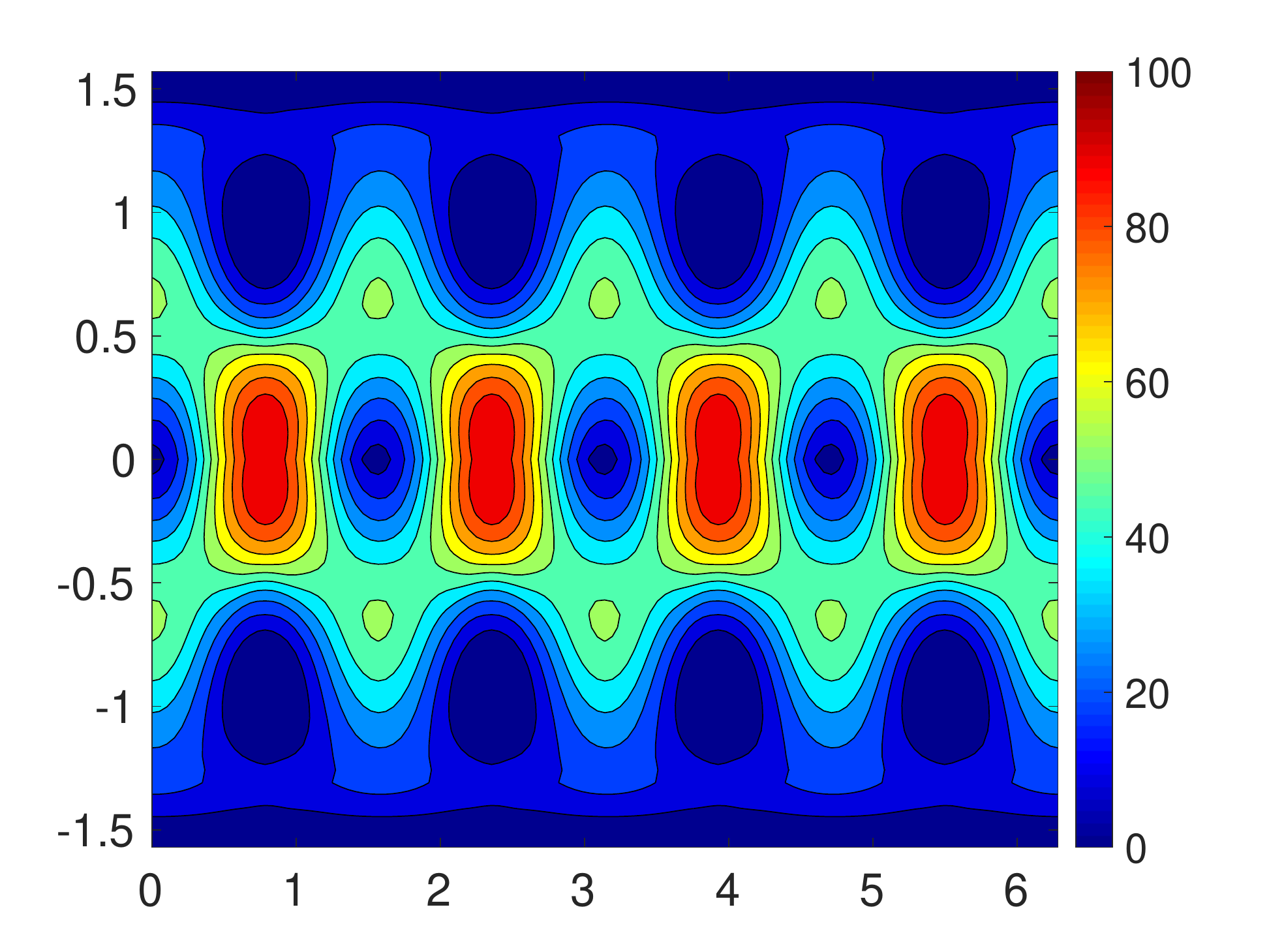}
		\subcaption{$u$ at $t=8$ days.}
	\end{subfigure}
	\begin{subfigure}[t]{0.32\textwidth}
		\centering
		\includegraphics[width=\linewidth]{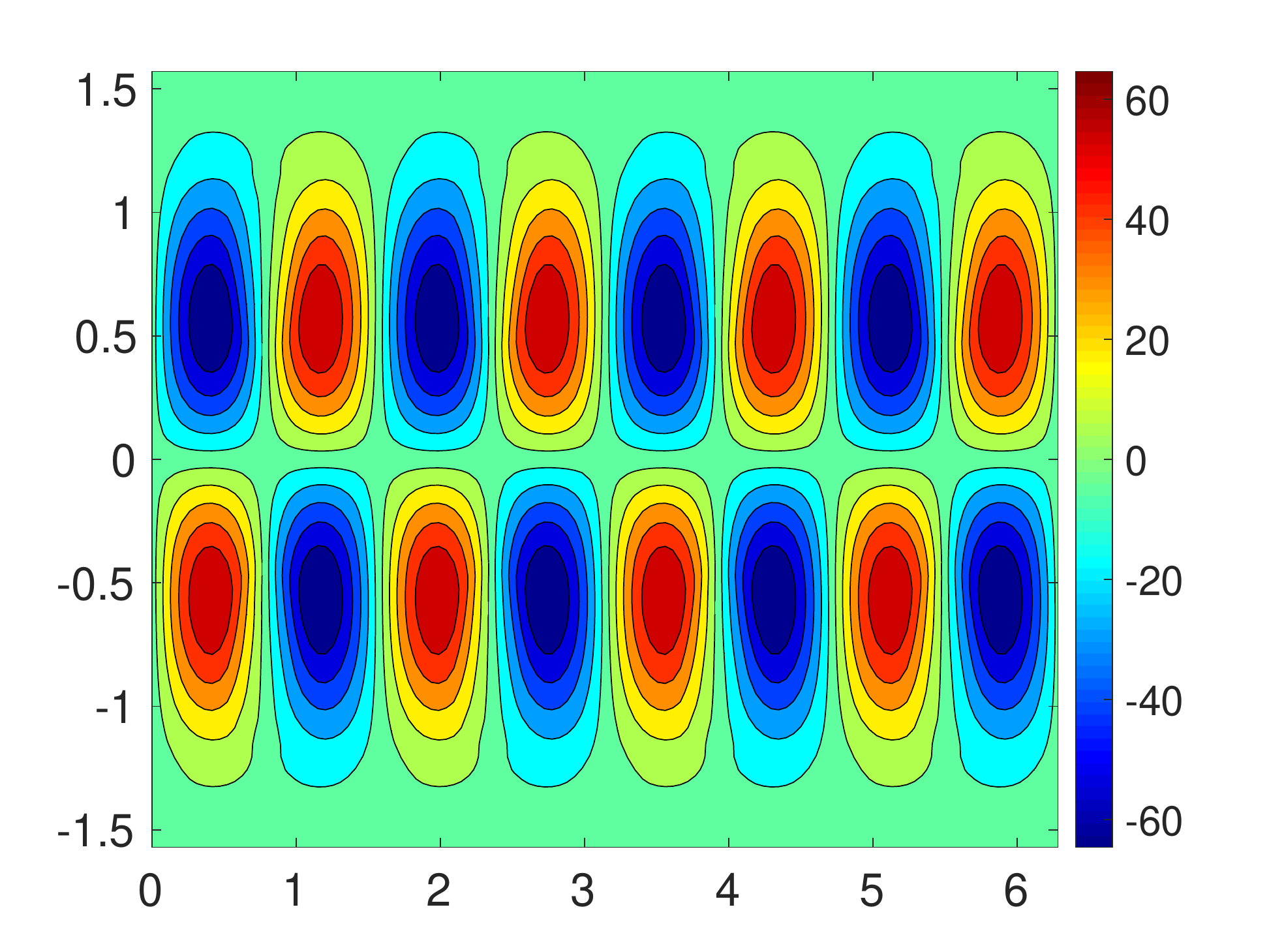}
		\subcaption{$v$ at $t=8$ days.}
\end{subfigure}

	\caption{8-day simulation of the Rossby-Haurwitz wave.}
	\label{fig:rossby-waves}
\end{figure}

The Rossby-Haurwitz wave (denoted as test case 6 in \cite{williamsonStandardTestSet1992}) consists of a large-scale planetary wave that mimics the high/low-pressure systems typical of mid-latitude weather patterns.
The test case considers initial data that would result in a stable solution for the barotropic vorticity equation, evolving from west to east without changing shape. It is known that this configuration is ultimately unstable --- see, for example, the discussion in \cite{thuburn2000numerical} --- but this instability only arises on a relatively long time scale. Therefore, it is customary to
assess the quality of numerical methods based on their capability to reproduce a stable eastward moving pattern for several days. In Figure \ref{fig:rossby-waves}, we see the results of an 8-day simulation of the Rossby-Haurwitz wave on a 40x20 grid using the RK4 method in time with $\Delta t=4$ s and $p=3$ in space.

    It can be observed that
    the numerical solution indeed evolves from west to east while maintaining a close resemblance with the initial shape and that the simulated pattern is in good agreement with reference solutions, see e.g., \cite{tumolo:2015}.
    
\section{Performance} \label{sec:performance}

In this section, we present performance benchmarks of the G4GT and GT4Py implementations for the SWE in both planar and spherical geometry. In order to provide an adequate comparison, we break down the G4GT implementation into three computing blocks which allow us to better analyze the complexity and compare the execution time of the temporal loop for the various backends supported by GT4Py. In both cases, the time spent in the precomputation steps is neglected, as it becomes negligible for long simulation periods.

The G4GT simulations have been run on the compute nodes of {\it Piz Daint} at CSCS, using an Intel\textsuperscript{\textregistered} Xeon\textsuperscript{\textregistered} E5-2690 v3, 12-core processor (single node), characterized by a peak memory bandwidth of 68 GB/s.
On the other hand, the GT4Py benchmarks were performed on a different partition of {\it Piz Daint} with the CPU code executed on two 18-core Intel\textsuperscript{\textregistered} Xeon\textsuperscript{\textregistered} E5-2695 v4 @ 2.10GHz processor (each with 77 GB/s peak memory bandwidth) and the GPU code on an NVIDIA\textsuperscript{\textregistered} Tesla\textsuperscript{\textregistered} P100 with 16GB of memory (540 GB/s peak memory bandwidth).  

\subsection{G4GT performance evaluation} 
 
The geostrophic adjustment setup presented in the previous section can also be used to evaluate the performance of the method and G4GT in general. Unless stated otherwise, the physical and numerical parameters are therefore kept unchanged, including a grid consisting of $50$x$50$ elements and a time step of $100 \ \rm s$.

\begin{figure}[htbp]
	\centering
	\begin{subfigure}{0.48\textwidth}
		\centering
		\includegraphics[width=\linewidth]{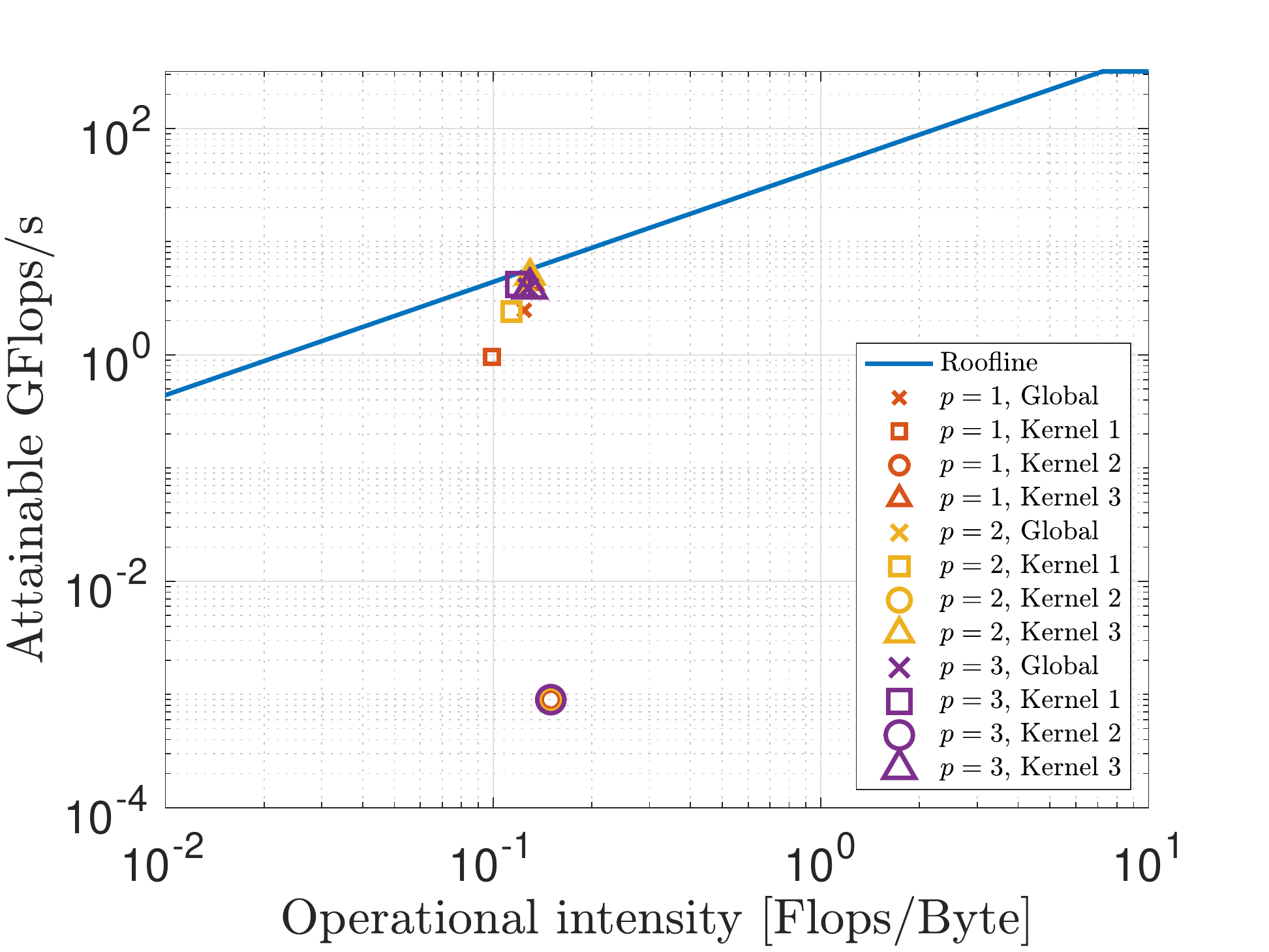}
		\subcaption{Global range}
		\label{fig:SWE_CPU_nozoom}
	\end{subfigure}
	\begin{subfigure}{0.48\textwidth}
		\centering
		\includegraphics[width=\linewidth]{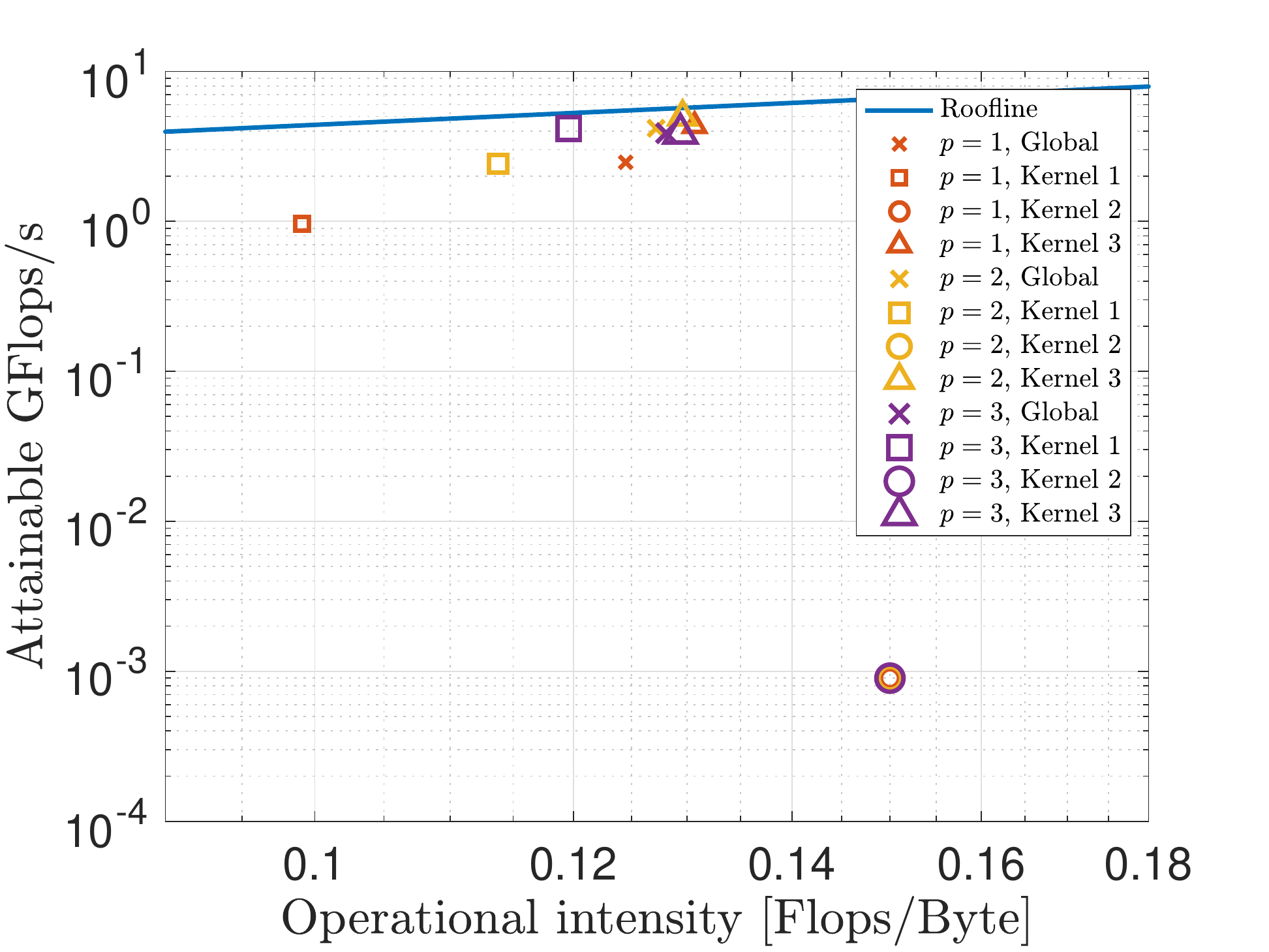}
		\subcaption{Zoom into the interval $[0.09,0.18]$}
		\label{fig:SWE_CPU_zoom}
	\end{subfigure}
	\caption{Performance evaluation of the G4GT implementation of the SWEs in a Cartesian geometry with different spatial degrees $p$, RK4 scheme.}
	\label{fig:SWE_CPU_RK4}
\end{figure}

The performance evaluation is done using the Roofline model \cite{Williams2009}. This is based on the \textit{operational intensity}, i.e., the number of floating-point operations (flops) per byte of DRAM traffic, and the \textit{attainable Gflops per second}, i.e., the concrete performance measure. Here, the DRAM traffic takes into account the bytes that are read from/written in the main memory after the filter of the cache hierarchy. Because of hardware limits, the attainable flops per second cannot go beyond a fixed threshold, determined by the peak memory bandwidth and the peak floating point performance. In practice, the actual threshold is determined by running benchmark cases, such as the (bandwidth-limited) \textit{stream} or the (computationally-limited) \textit{linpack} benchmark. In our case, these give a memory bandwith upper bound of $44 \ GB/s$ and a peak performance limit of $318 \ GFlops/s$ \cite{Calore2020}.
Thus, for a given operational intensity, an efficient implementation in terms of performance should attain values close to the determined limit. In our analysis, we decided to ignore the cache effects. In other words, every access to a variable is considered for the computation of the required bytes. This is in contrast with the definition provided by the model, but a precise estimate of the DRAM traffic is far from an easy task and, in the G4GT framework, no tool is available to appopriately measure it.
The matrix-vector multiplication, which is the central operation in the DG implementation, is bandwidth-limited and achieves a performance \cite{discacciati2018implementation} somewhat below the leftmost (rising) roofline. \\
Based on the way in which the code is structured \cite{discacciati2018implementation}, we can recognize three different kernels:
\begin{enumerate}
	\item {\bf Common part:} The nodal values for the solution and the flux function are computed.
	\item {\bf Rusanov fluxes:} The boundary fluxes are computed, and the boundary conditions are applied. This requires communication among neighboring elements.
	\item {\bf Main computation:} The right-hand side is assembled, and the solution is updated.
\end{enumerate}
The results for varying polynomial degree $p$ are reported in Figure \ref{fig:SWE_CPU_RK4}, which compares the performances of the global program and the kernels separately. As a complement to Figure \ref{fig:SWE_CPU_RK4}, Table \ref{table:SWE_times_CPU} breaks down the computational times.
Looking at the overall performance, we observe that no significant variations in the operational intensities are present. This is because variations in the polynomial order lead to similar changes in the number of floating point operations and memory traffic. 
However, it appears that performances obtained with linear basis functions are slightly lower than higher-order polynomials. The total number of operations might not be large enough to attain the expected asymptotic values, causing deviations from the optimal performance. 

Looking at the kernels independently, for high values of $p$ the third kernel has the most significant influence on the overall performance. This is expected since it includes the majority of the computations. Specifically, the assembly of the internal integral is the most intensive part, both in terms of resources and time. 
On the other hand, the second kernel always has a low computational cost. This is not surprising, as only boundary quantities are involved. The performance of this kernel in terms of floating point operations per second are consistently low and do not vary with $p$. Since it is the only phase that involves exchanges between neighboring elements, it is reasonable that cache misses or inefficient memory accesses are present. 
No particular trend is observed for the first kernel, except for $p=1$, in which this kernel has the dominant effect on global performance.
\begin{table}[htbp]
	\centering
	\begin{tabular}{lcccc}\toprule
		$p$ & Global [$\rm s$] & Kernel 1 [$\rm s$] & Kernel 2 [$\rm s$] & Kernel  3 [$\rm s$]\\\midrule
		1 & 15.17 & 7.42 & 1.45 & 6.31 \\
		2 & 78.98 & 21.22 & 3.10 & 54.66 \\
		3 & 430.21 & 54.97 & 5.81 & 369.43 \\\bottomrule
	\end{tabular}
	\caption{Computational times of the G4GT implementation of the SWEs in a Cartesian geometry with different spatial degrees $p$, RK4 scheme.}
	\label{table:SWE_times_CPU}
\end{table}

\subsection{GT4Py performance evaluation}

For the GT4Py implementation, we consider both performance scalability while increasing the horizontal problem size 
as well as increasing the number vertical layers while holding the horizontal size constant.

\subsubsection{Horizontal scaling}
\begin{figure}[htb]
        \centering
        \includegraphics[width=0.5\textwidth]{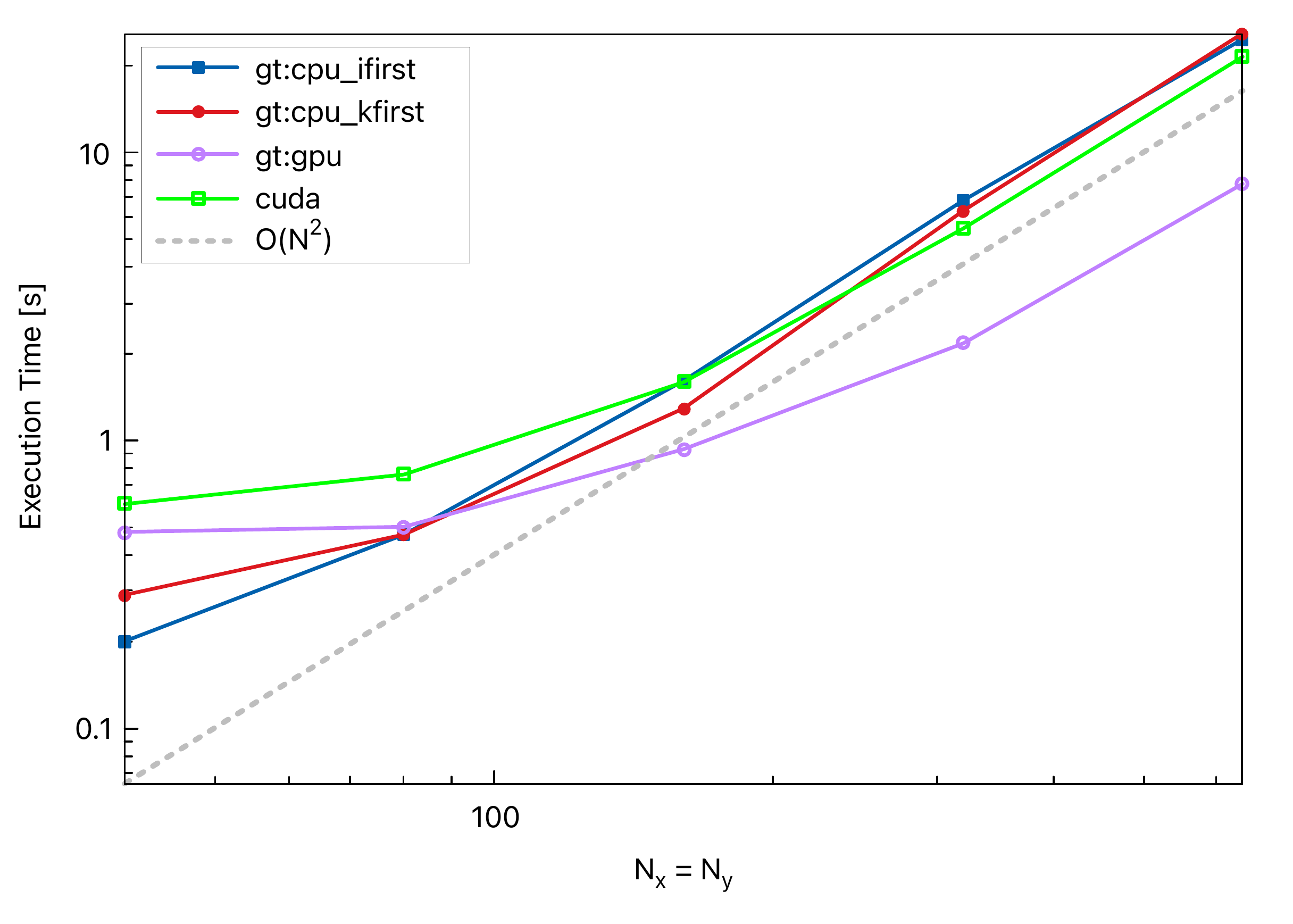}
        \caption{Benchmark of the execution time of the GT4Py backends  with increasing problem size.
        Each subsequent data point doubles the grid points in $x$ and $y$ and thus quadruples the total number of grid points.}
        \label{fig:nz_scaling}
    \end{figure}

    We study the scaling of the runtime of the application with increasing horizontal resolution. In the benchmark, we use a 4th-order scheme in space which corresponds to a vector of size 16 stored at each grid point (\texttt{data\_dims = 16}).
    Figure~\ref{fig:nz_scaling} compares the runtimes of the various backends. The \texttt{gt:cpu\_ifirst} backend was used as the baseline, since this is the best performing CPU backend.
    
    Each successive data point doubles the number of grid points in both horizontal directions.
    Thus, the expected asymptotic scaling is quadratic since doubling the number of grid points in both horizontal directions results in a four-fold increase in the total number of grid points.

    In this experiment, the two CPU backends, powered by the GridTools framework, have virtually identical execution times.
    For small problem sizes, the two GPU backends (namely the \texttt{cuda} and \texttt{gt:gpu}) perform worse than the CPU backends.
    This is due to the under-utilization of the GPU resources for small-scale problems which are not able to fully saturate the device.
    This is illustrated in the delayed asymptotic scaling of the GPU code when compared to the CPU resulting in better performance for larger problems.
    When comparing the two GPU backends, we observe that the optimizations provided by the GridTools framework (included in the \texttt{gt:gpu} backend) yields significant better performing code than the naive CUDA backend.
    Note that we were limited to presenting a maximum problem size of 640x640 due to memory constraints on the GPU.

    Table \ref{table:speedup} summarizes the speedup observed versus the \texttt{gt:cpu\_ifirst} baseline on the largest problem size.
    \begin{table}[htb]
        \centering
        \begin{tabular}{ccccc}\toprule
             & \texttt{gt:cpu\_kfirst} & \texttt{gt:cpu\_ifirst} & \texttt{cuda} & \texttt{gt:gpu} \\\midrule
            Speedup Factor & 0.952 & 1.00 & 1.15 & 3.19
        \end{tabular}
        \caption{Speedup of GT4Py backends vs reference \texttt{gt:cpu\_ifirst} implementation on 640x640 grid.}
        \label{table:speedup}
    \end{table}
One observation from Table \ref{table:speedup} is that the fastest GPU backend results in a speedup factor of $\sim 3.2$ versus the fastest CPU backend. Since we know from Figure~\ref{fig:SWE_CPU_nozoom} that the performance is limited by memory bandwidth, we expect the speedup to mirror the ratio of bandwidths listed in Section \ref{sec:performance} for the P100 GPU to two Intel Broadwell processors, namely, 560 GB/s : 2x77 GB/s $\simeq 3.6$.

    
    
\subsubsection{Vertical scaling}

     In this Section we try to assess the potential performance of GT4Py on a 3-dimensional problem by considering a set of
    decoupled 2-dimensional SWE problem  copied in the vertical direction and solved in parallel.  This configuration increases the computational load and resembles to some extent those of low order finite difference/finite volume discretizations of 3-dimensional problems in atmospheric modeling. However, it is substantially different from a full 3-dimensional DG discretization, since all the local matrices that arise correspond to 2-dimensional rather than 3-dimensional elements. The resulting algorithm scales linearly with the number of vertical levels.
    Considering that all levels can be solved independently, this problem is embarrassingly parallel, and we might expect performance benefits on the GPU compared to the CPU.

\begin{figure}[htb]
        \centering
        \includegraphics[width=0.5\textwidth]{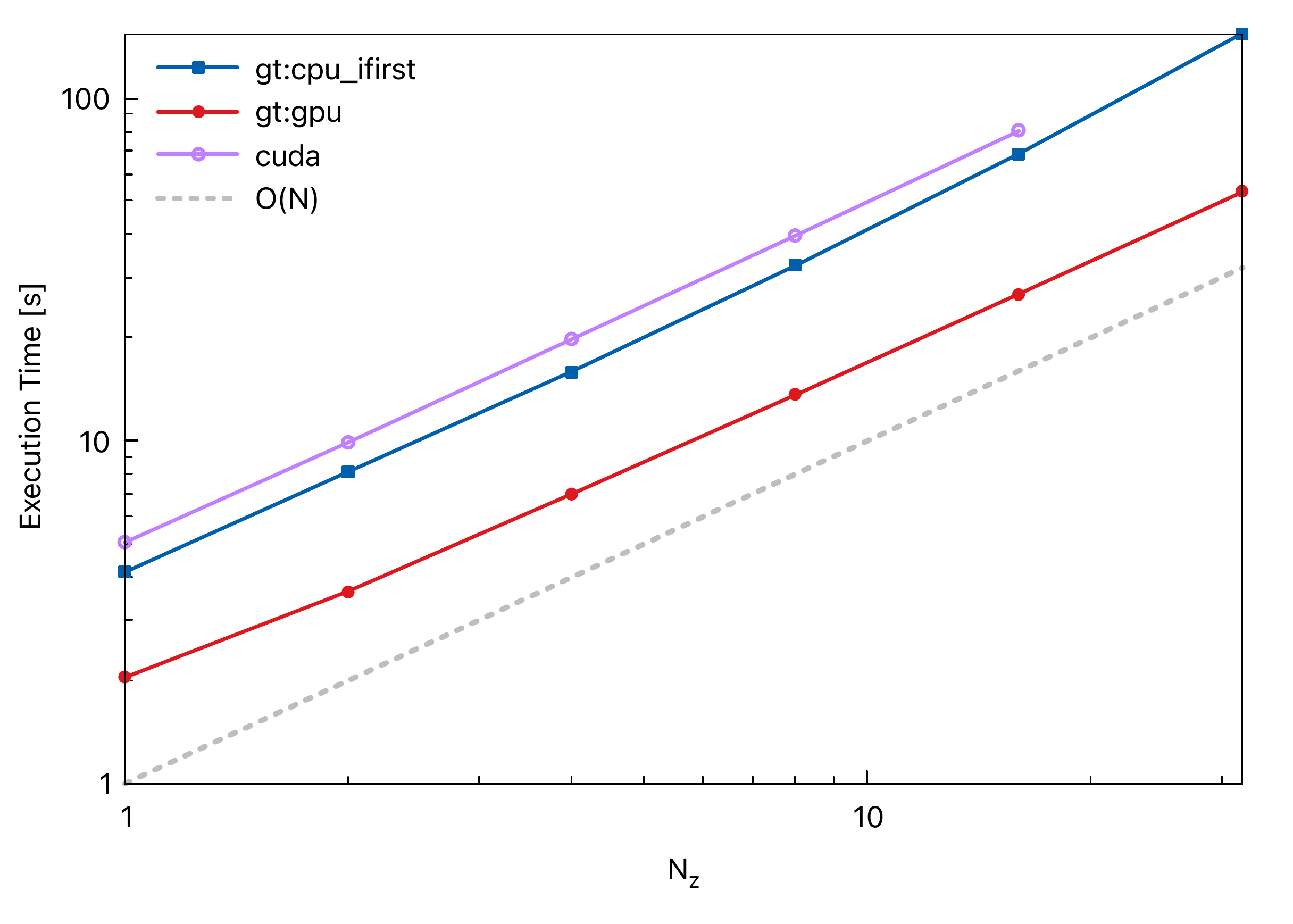}
        \caption{Benchmark of best-performing CPU and GPU backends in addition to the CUDA backend.
        The plot depicts execution time with respect to the number of identical vertical problems solved in parallel.
        The final data point for the CUDA backend is unavailable due to the memory limit reached on GPU.}
        \label{fig:nx300}
    \end{figure}

    Figure \ref{fig:nx300} illustrates the execution time of a 4th-order DG scheme on a 300x300 grid with increasing vertical levels.
    Surprisingly, we do not observe any scaling benefits for this experiment on the GPU.
    Indeed, all backends exhibit asymptotic scaling from the first data point, indicating that the hardware's resources are fully saturated.
    This is most likely due to the 2-dimensional problem solved in each level being too large and fully occupying the memory bandwidth of the GPU.
    We observe that the CUDA backend performs even worse than the GridTools CPU backend.
    Moreover, it suffers from poor memory utilization compared to the GridTools GPU implementation, as the last data point could not be gathered due to the memory capacity of the GPU being reached.
    

\subsubsection{GPU profiling}
    In this section, we present the results of a brief profiling that was carried out for the GridTools and CUDA GPU backends using the NSight\textsuperscript{\textregistered} code analysis tools from NVIDIA\textsuperscript{\textregistered}.

\begin{figure}[htb]
        \centering
        \includegraphics[width=0.7\textwidth]{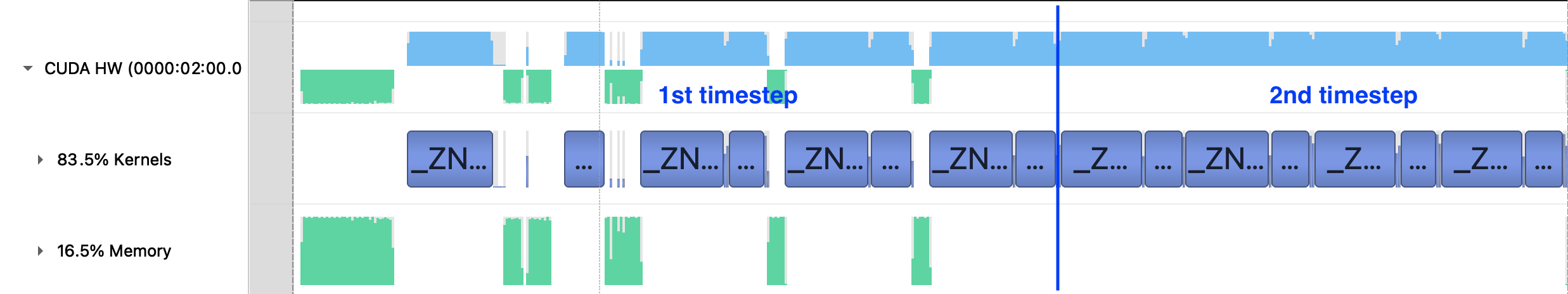}
        \caption{Extract of NSight Systems profile of \texttt{gt:gpu} backend.
        The profile indicates consistent executions of kernels on GPU and no communication between the host and device after the first time step.}
        \label{fig:profile_timeline}
    \end{figure}

    We profile the 4th-order DG scheme for the 2-dimensional SWE on a 640x640 grid.
    Figure \ref{fig:profile_timeline} illustrates an extract of the profile timeline for the first two time steps of the simulation.
    The green sections indicate memory transfers from the host to the device, while the blue areas highlight kernel executions.
    To achieve good performance on the GPU, communication between the device and the host needs to be reduced as much as possible.
    Considering the memory row in Figure \ref{fig:profile_timeline}, we only observe communication (more precisely host-to-device transfers) in the first time step.
    These are necessary to transfer the initial fields to the GPU.
    Subsequently, the GPU is able to execute kernels consistently without having to wait for unnecessary communication from the host.


\section{Conclusions}
\label{sec:conclusions}

We have presented two implementation examples of a high-order DG method in the framework of the G4GT and GT4Py domain-specific languages, respectively.
The G4GT and GT4Py implementations illustrate that a DSL designed for finite difference/volume methods on rectangular grids, with some extensions, can be reused to achieve the implementation of a DG solver.
The related performance analysis illustrated clearly that the DG method is limited by the memory bandwidth, similar to the sparse matrix-vector (SpMV) stencil discussed in  \cite{Williams2009}.

Despite being only a proof of concept, the G4GT DSL has several advantages for the end user over GT, with the main one being the simplicity in its use. Most of the back-end optimizations and the definition of appropriate data structures are hidden in the underlying GT implementation. Therefore, the user can exploit the full capability of GT in a broader range of discretization schemes without delving into the details required to achieve computational efficiency, in complete agreement with the requirements of a good DSL.
Several technical reasons, including the need for external libraries and the inability to construct nested functors, as well as the emergence of Python as a programming language, were responsible for the termination of the development of G4GT and the migration to GT4Py.

Despite the user-friendly nature of GT4Py and the increasing support and functionalities available, we had to expand the support for higher-dimensional fields by supplementing the frontend with an new, clean syntax. This allowed us to port an original Matlab implementation to GT4Py with relative ease. Moreover, we saw that we could automatically exploit available accelerators without having to modify the source code.
Nonetheless, due to limitations of the functionality of higher-dimensional fields, we were left with writing a lot of boilerplate code in GT4Py.
At the time of the implementation, slicing was not supported, which would have allowed us to condense the conserved variables into a large matrix instead of a series of separate vectors. In addition, there was no method  for calling functions on higher-dimensional fields which would eliminate the need to copy the same functionality to different stencils.

In conclusion, the new GridTools support for higher-dimensional fields opens the door for implementing advanced numerical schemes like DG in a DSL environment.
However, our front-end changes should be complemented with corresponding back-end optimizations to achieve better performance. 

\section*{Acknowledgments}

The Swiss National Supercomputing Centre (CSCS) funded the four-month internships of the first two authors, N.D. and K.S.. We would like to thank Linus Groner, Till Ehrengruber, Mauro Bianco and Christopher Bignamini of CSCS for their gracious support during both internships. L.B. was partially supported by the ESCAPE-2 project, European Union’s Horizon 2020
Research and Innovation Programme (Grant Agreement No. 800897).

\bibliographystyle{plainnat}
\bibliography{biblio}

\end{document}